\documentclass[A4paper]{article}

\newcommand{\footremember}[2]{%
    \footnote{#2}
    \newcounter{#1}
    \setcounter{#1}{\value{footnote}}%
}
\newcommand{\footrecall}[1]{%
    \footnotemark[\value{#1}]%
} 

\usepackage{style}
\begin{document}

\title{The Mixed Virtual Element Discretization for highly-anisotropic problems: the role of the boundary degrees of freedom}

\author{S. Berrone\footremember{trailer}{Department of Applied Mathematics, Politecnico di Torino, Italy
  (stefano.berrone@polito.it, stefano.scialo@polito.it, gioana.teora@polito.it).}, S. Scial\`o\footrecall{trailer}{}, G. Teora\footrecall{trailer}{}}
\maketitle

\begin{abstract}
In this paper, we discuss the accuracy and the robustness of the mixed Virtual Element Methods when dealing with highly-anisotropic diffusion problems. In particular, we analyze the performances of different approaches which are characterized by different sets of both boundary and internal degrees of freedom in presence of a strong anisotropy of the diffusion tensor with constant or variable coefficients. A new definition of the boundary degrees of freedom is also proposed and tested. 
\end{abstract}


\section{Introduction}

The Virtual Element Method (in short VEM) \cite{LBe13, secondVEM} is a generalization of the Finite Element
Method (FEM in short) that can easily handle general polytopal meshes and high-order methods. The major difference with the FEM is that the VEM space contains suitable non-polynomial functions. For this reason, the standard VEM discrete bilinear form is the sum of a consistency part ensuring accuracy and of a stabilization term enforcing the coercivity. In particular, the choice of the stabilization term remains a critical part of the VEM construction \cite{russo2023quantitative, berrone2023lowest} and it is usually problem-driven. Furthermore, the stabilization term may have possible negative effects on the conditioning of the system \cite{mascotto,DASSI20183379} and may become an issue in highly anisotropic diffusion problems due to its isotropic nature.
Moreover, we recall that in order to build high-order methods, it is crucial to employ a well-conditioned polynomial basis in the definition of the internal degrees of freedom in order to obtain reliable solutions. Indeed, the advantages of using $L^2$-orthogonal polynomial bases against the standard monomial one have largely been proved both for the primal version of the method \cite{mgs_basis,Sbe17,mascotto,DASSI20183379, berrone2023improving} and for its mixed formulation \cite{berrone2023orthogonal}.

In this paper, we want to test the accuracy and the robustness of the mixed Virtual Element Method when dealing with highly anisotropic diffusion tensors. For this purpose, we propose different kinds of degrees of freedom and test them against different choices of the stabilization term for a set of benchmark anisotropic diffusion problems. In particular, we introduce a new set of boundary degrees of freedom which are defined as moments up to degree $k \geq 0$ against an $L^2([0,1])$-orthonormal polynomial basis in order to analyze the role of the boundary degrees of freedom in the conditioning and in the accuracy of the methods. Numerical experiments show that this choice of boundary degrees of freedom generally leads to a downward shift of the error curves. However, this approach does not result in an improvement of the condition number of the system matrix in all the test cases.

The outline of the paper is as follows. In Section \ref{sec:modelproblem} we present the model problem. In Section \ref{sec:localspace}, after introducing the local mixed virtual element spaces and different sets of the local degrees of freedom, we define the mixed VE formulation of the problem. In Section \ref{sec:stabilization}, we describe the main properties and discuss possible choices for the stabilization term. Finally, in Section \ref{sec:experiments} we test all the proposed approaches through different benchmark problems which are characterized by highly anisotropic diffusion tensors, with both constant and variable coefficients.

\section{The model problem}\label{sec:modelproblem}
Let $\Omega \subset \R^2$ be a bounded convex polytopal domain with boundary $\Gamma$ and let $\nn_{\Gamma}$ be the outward unit normal vector to the boundary. 
Let us consider a tensor $\D(\x) \in \R^{2\times 2}$ which is bounded, measurable, symmetric and strongly elliptic on $\Omega$, i.e. there exist $\D_{min},\ \D_{max}$, independent on $\vv$ and $\x$, such that
\begin{linenomath}
\begin{equation*}
\D_{min} \norm{\vv(\x)}^2 \leq \vv(\x) \cdot \D(\x) \vv(\x) \leq \D_{max} \norm{ \vv(\x) }^2,
\end{equation*}
\end{linenomath}
holds for every $\vv \in H_{0,\Gamma_N}(\div ;\Omega) = \{v \in H(\div; \Omega): \vv \cdot \nn_{\Gamma_N} = 0\}$ and for almost every $\x \in \Omega$,
where $\norm{}$ denotes the euclidean norm.
Given $f\in L^2(\Omega)$, $g_D \in H^{1/2}(\Gamma_D)$ and $g_N \in L^2(\Gamma_N)$, we consider the following diffusion problem
\begin{equation}
\begin{cases}
\div \left(- \D \nabla p\right) = f & \text{in } \Omega\\
p = g_D & \text{on } \Gamma_D\\
- (\D \nabla p) \cdot \nn_{\Gamma_N} = g_N & \text{on } \Gamma_N
\end{cases},
\label{eq:primalForm}
\end{equation}
where $\Gamma_D$ and $\Gamma_N$ such that $\Gamma_D \cup \Gamma_N = \Gamma$ and $\vert \Gamma_n \cap \Gamma_N \vert = 0$ denote the Dirichlet and the Neumann boundary, respectively.
In particular, in the following, we focus on diffusion problems with a diffusion tensor of the form
\begin{equation}
\D(\x) = \bm{R}(\x) \begin{bmatrix}
\D_{max} & 0 \\
0 & \D_{min}
\end{bmatrix} (\bm{R}(\x))^T,
\end{equation}
which is characterized by a high anisotropic ratio, i.e. the ratio between the smallest and largest eigenvalues of the diffusion tensor.

Introducing the velocity space $\V = H_{0,\Gamma_N}(\div ;\Omega)$ and the pressure space $\Q=L^2(\Omega)$, the mixed variational formulation of \eqref{eq:primalForm} reads:
\begin{equation}
\begin{cases}
\text{Find } (\uu_0,p)\in \V \times \Q \text{ such that } \uu = \uu_0 + \uu_N \text{ and } p \text{ satisfy} & \\
\scal[\Omega]{\D^{-1}\uu}{\vv} - \scal[\Omega]{p}{\div \vv} = -\langle g_D, \vv \cdot \nn_{\Gamma_D} \rangle_{\pm\frac{1}{2}, \Gamma_D}  & \forall \vv \in \V\\
\scal[\Omega]{\div \uu}{q} = \scal[\Omega]{f}{q} & \forall q\in \Q
\end{cases},
\label{eq:mixedForm}
\end{equation}
where $\uu_N \in H(\div;\Omega)$ is a chosen function that satisfies $\uu_N \cdot \nn_{\Gamma_N} = g_N$ and $\langle \cdot, \cdot \rangle_{\pm\frac{1}{2}, \Gamma_D}$ denotes the duality paring between $H^{-1/2}(\Gamma_N)$ and $H^{1/2}(\Gamma_N)$. 

\section{The mixed Virtual Element Space}\label{sec:localspace}

Now, let us consider a decomposition $\Th$ of $\Omega$ in star-shaped polygons $E$, where $h$, as usual, is set to be the maximum diameter of elements $E \in \Th$. We further denote by $\Eh[E]$ the set of edges of an element $E \in \Th$.

For any integer $k \geq 0$, we define the local virtual element space related to the velocity variable $\uu$ as 

\begin{multline*}
\Vh[E]{k} = \Big\{ \vv \in H(\div; E) \cap H(\rot ; E):\ \vv \cdot \nn_e \in \Poly{k}{e}\forall e \in \Eh[E],\\ \div \vv \in \Poly{k}{E},\  \rot \vv \in \Poly{k-1}{E}\Big\},
\end{multline*}
and the local virtual element space related to the pressure variable $p$ as $\Qh[E]{k} = \Poly{k}{E}$, which is the space of the polynomials of order up to $k$ on $E$ \cite{Mixed}.

The choice of the degrees of freedom in the local pressure space $\Qh[E]{k}$ is trivial: the degrees of freedom of a function $p \in \Qh[E]{k}$ are its coefficients with respect to the polynomial basis chosen as the basis for $ \Poly{k}{E}$.
The standard polynomial basis for $\Poly{k}{E}$ used in the VEM construction \cite{MixedImplem} is given by the set of the $n_{k} = \dim \Poly{k}{E} = \frac{(k+1)(k+2)}{2}$ bi-dimensional scaled monomials, i.e.
\begin{equation}
\M{k}{E} = \Big\{ m_{\alpha} = \left(\frac{\x-\x_E}{h_E}\right)^{\bm{\alpha}}: \alpha= \ell(\bm{\alpha})\ \forall \alpha = 1,\dots,n_{k} \Big\}
\end{equation}
where $\x_E$ and $h_E$ are the centroid and the diameter of the polygon $E$, respectively, and $\ell$ is the function $\N^2 \to \N$ which maps
\begin{linenomath}
\begin{equation*}
(0,0) \mapsto 1,\quad (1,0) \mapsto 2,\quad (0,1) \mapsto 3,\quad (2,0) \mapsto 4,\cdots
\end{equation*}
\end{linenomath}
A more robust choice is represented by the set of the $L^2(E)$-orthonormal polynomials $\QPoly{k}{E} = \{q_{\alpha}\}_{\alpha=1}^{n_k}$ introduced in \cite{mgs_basis,Sbe17,mascotto} for the primal version of the method and then tested in the mixed case in \cite{berrone2023orthogonal}. This orthonormal polynomial basis is defined as
\begin{equation}
 q_{\beta} = \sum_{\gamma=1}^{n_{k}} \mathbf{L}^k_{\beta \gamma}  m_{\gamma},\ \forall \beta=1,\dots,n_{k},
 \label{eq:qbasis}
\end{equation}
where $\mathbf{L}^k \in \R^{n_k \times n_k}$ is built by applying twice the modified Gram Schmidt algorithm to the monomial Vandermonde matrix related to a proper quadrature formula on $E$.

\subsection{The Degrees of Freedom for the velocity variable}

In order to define the local degrees of freedom for the local velocity space $\Vh[E]{k}$, we need to introduce the following polynomial spaces.
We introduce the (vector) polynomial space
\begin{equation}
\GPoly{k}{\nabla,m}{E} = \nabla \M{k+1}{E} = \Big\{\g^{\nabla,m}_{\alpha}\Big\}_{\alpha=1}^{n^{\nabla}_k} \subset \VPoly{k}{E},
\end{equation}
and the set $\GPoly{k}{\perp,m}{E} = \Big\{\g^{\perp,m}_{\alpha}\Big\}_{\alpha=1}^{n^{\perp}_k}$ which is defined in such a way
\begin{linenomath}
\begin{equation*}
\VPoly{k}{E} = \GPoly{k}{\nabla,m}{E} \oplus \GPoly{k}{\perp,m}{E},
\end{equation*}
\end{linenomath}
with $n^{\nabla}_k = n_k + (k+1)$ and $n^{\perp}_k = n_{k}- (k+1)$.
The set $ \GPoly{k}{m}{E} = \GPoly{k}{\nabla,m}{E} \cup \GPoly{k}{\perp,m}{E}$ represents a (vector) polynomial basis for $\VPoly{k}{E}$ which allows to easily define the set of local degrees of freedom in the mixed VEM framework \cite{MixedImplem}.
Let us denote by $\GPoly{k}{\nabla,\overline{q}}{E} = \{\g^{\nabla,\overline{q}}_{\alpha}\}_{\alpha=1}^{n^{\nabla}_{k}} \subset \VPoly{k}{E}$ the set of (vector) polynomials
\begin{equation}
 \g^{\nabla,\overline{q}}_{\alpha} = \sum_{\beta=1}^{n_k^{\nabla}} \mathbf{L}^{\nabla,k}_{\alpha \beta} \nabla q_{\beta+1},\ \forall \alpha=1,\dots,n^{\nabla}_k
 \label{eq:gnablaortho}
\end{equation}
such that
\begin{linenomath}
\begin{equation*}
\scal[E]{\g^{\nabla,\overline{q}}_{\alpha}}{\g^{\nabla,\overline{q}}_{\beta}} = \delta_{\alpha \beta},\quad \forall \alpha,\beta =1,\dots,n^{\nabla}_{k},
\end{equation*}
\end{linenomath}
which is obtained by orthonormalizing the gradients of polynomials belonging to $\QPoly{k+1}{E}$ throughout the modified Gram-Schmidt algorithm.
Now, we define $\GPoly{k}{\perp,\overline{q}}{E} = \{\g^{\perp,\overline{q}}_{\alpha}\}_{\alpha=1}^{n^{\perp}_{k}}$ as the $L^2(E)$-orthogonal complement of $\GPoly{k}{\nabla,\overline{q}}{E}$ in $\VPoly{k}{E}$, which is chosen such that
\begin{linenomath}
\begin{equation*}
\scal[E]{\g^{\perp,\overline{q}}_{\alpha}}{\g^{\perp,\overline{q}}_{\beta}} = \delta_{\alpha \beta},\quad \forall \alpha,\beta =1,\dots,n^{\perp}_{k}.
\end{equation*}
\end{linenomath}
Further details about the construction of this basis can be found in \cite{berrone2023orthogonal}. Here, it was shown that it is advisable to choose the set
\begin{equation}
\GPoly{k}{\overline{q}}{E}= \GPoly{k}{\nabla,\overline{q}}{E} \cup \GPoly{k}{\perp,\overline{q}}{E},
\end{equation}
as the (vector) polynomial basis for $\VPoly{k}{E}$ in order to reduce the ill-conditioning of the system matrix and to obtain more accurate and reliable solutions for high values of the local polynomial degree and in presence of badly-shaped polygons.

Now, let us introduce a quadrature formula $\mathbb{S}^Q = \{(s_j^Q,w_j^Q)\}_{j=1}^{N^Q}$ of order $2(k+1)$ with $N^Q \geq k+2$ nodes on the interval $[0,1]$. 
We define the one-dimensional $L^2([0,1])$-orthonormal polynomial basis $\QPoly{k+1}{[0,1]} = \{t_1,\dots,t_{k+1},t_{k+2}\}$ for $\Poly{k+1}{[0,1]}$ by applying the modified Gram-Schmidt algorithm with reorthogonalization to the Vandermonde matrix $\mathbf{V}^{\mathbb{S}^Q} \in \R^{N^Q \times (k+2)}$ related to the one-dimensional monomial basis $\{1,s,\dots,s^k,s^{k+1}\}$ and the quadrature formula $\mathbb{S}^Q$. More precisely, we perform sequentially
\begin{linenomath}
\begin{equation*}
\mathbf{V}^{\mathbb{S}^Q} = \mathbf{Q}^{\mathbb{S}^Q}_1 \mathbf{R}^{\mathbb{S}^Q}_1, \quad \mathbf{R}^{\mathbb{S}^Q}_1  \in \R^{(k+2)\times (k+2)},\ \mathbf{Q}^{\mathbb{S}^Q}_1  \in \R^{N^Q \times (k+2)}: (\mathbf{Q}^{\mathbb{S}^Q}_1)^T \mathbf{Q}^{\mathbb{S}^Q}_1 = I
\end{equation*}
\begin{equation*}
\sqrt{\mathbf{W}^{\mathbb{S}^Q}} \mathbf{Q}_1^{\mathbb{S}^Q} = \mathbf{Q}^{\mathbb{S}^Q}_2 \mathbf{R}^{\mathbb{S}^Q}_2, \quad \mathbf{R}^{\mathbb{S}^Q}_2  \in \R^{(k+2)\times (k+2)},\ \mathbf{Q}^{\mathbb{S}^Q}_2  \in \R^{N^Q \times (k+2)}: (\mathbf{Q}^{\mathbb{S}^Q}_2)^T \mathbf{Q}^{\mathbb{S}^Q}_2 = I,
\end{equation*}
\end{linenomath}
where $\mathbf{W}^{\mathbb{S}^Q} \in \R^{N^Q \times N^Q}$ is the diagonal matrix of quadrature weights, and then we define
\begin{equation}
t_j = \sum_{i=1}^{k+2} \mathbf{L}^{\mathbb{S}^Q,k+1}_{ji} s^i,\quad \forall j=1,\dots,k+2,
\end{equation}
where $\mathbf{L}^{\mathbb{S}^Q} = (\mathbf{R}^{\mathbb{S}^Q}_2 \mathbf{R}^{\mathbb{S}^Q}_1)^{-T}$.

We remark that each polynomial in $\Poly{k+1}{e}$, $e\in\Eh[E]$, can be written in terms of polynomials in $\QPoly{k+1}{[0,1]}$ through an affine mapping $F:[0,1] \to e$.
Furthermore, we recall that the modified Gram-Schmidt algorithm is a hierarchical procedure, which means, for example,
\begin{linenomath}
\begin{equation*}
\QPoly{k}{[0,1]} = \{t_1,\dots,t_{k+1}\} \subset \QPoly{k+1}{[0,1]},
\end{equation*}
\end{linenomath}
is a basis for $\Poly{k}{[0,1]}$.

In $\Vh[E]{k}$, we define the set of local Degrees of Freedom (DOFs in short) as the union of 
\begin{enumerate}[label=\textbf{\arabic*.}]
\item the set of the boundary degrees of freedom which can be chosen as
\begin{enumerate}[label*=\textbf{\alph*)},ref=\textbf{1.\alph*)}]
\item \label{DOF:1a} the values of $\vv_h \cdot \nn_e$ in the $k+1$ Gauss quadrature points $\x_i^{e,Q}$ internal on each edge $e \in \Eh[E]$,
\end{enumerate}
or
\begin{enumerate}[label*=\textbf{\alph*)},ref=\textbf{1.\alph*)}]
\setcounter{enumii}{1}
\item \label{DOF:1b} the $k+1$ moments on each edge $e \in \Eh[E]$:
\begin{equation}
\int_0^1 \widehat{\vv_h \cdot \nn_e} t_j \vert e \vert,\quad \forall j=1,\dots,k+1,
\end{equation}
where $\vert e \vert$ represents the length of the edge $e$, while $(\widehat{\vv_h \cdot \nn_e})(s) = (\vv_h \cdot \nn_e)(F(s))$.
\end{enumerate}
\item the set of the internal degrees of freedom which can be chosen as the internal moments computed against
\begin{enumerate}[label*=\textbf{\roman*)},ref=\textbf{2.\roman*)}]
\item \label{DOF:2a}the sets of functions $\GPoly{k-1}{\nabla,m}{E}$ and $\GPoly{k}{\perp,m}{E}$:
\begin{equation}
\frac{1}{\vert E \vert} \int_E \vv_h \cdot \g^{\nabla,m}_{\alpha},\quad \forall \alpha = 1,\dots,n_{k-1}^{\nabla},
\label{eq:idof_nabla_mon}
\end{equation}
\begin{equation}
\frac{1}{\vert E \vert} \int_E \vv_h \cdot \g^{\perp,m}_{\alpha},\quad \forall \alpha = 1,\dots,n_k^{\perp},
\label{eq:idof_perp_mon}
\end{equation}
\end{enumerate}
or 
\begin{enumerate}[label*=\textbf{\roman*)},ref=\textbf{2.\roman*)}]\setcounter{enumii}{1}
\item \label{DOF:2b}the sets of functions $\GPoly{k-1}{\nabla,\overline{q}}{E}$ and $\GPoly{k}{\perp,\overline{q}}{E}$:
\begin{equation}
\frac{1}{\vert E \vert} \int_E \vv_h \cdot \g^{\nabla,\overline{q}}_{\alpha},\quad \forall \alpha = 1,\dots,n_{k-1}^{\nabla},
\end{equation}
\begin{equation}
\frac{1}{\vert E \vert} \int_E \vv_h \cdot \g^{\perp,\overline{q}}_{\alpha},\quad \forall \alpha = 1,\dots,n_k^{\perp},
\end{equation}
\end{enumerate}
where $\vert E \vert$ is the area of the polygon $E$.
\end{enumerate}

Let us denote by $\Ndof[E] = \dim \Vh[E]{k} = \#\Eh[E] (k+1)+ n^{\nabla}_{k-1} + n^{\perp}_k$ and let us introduce the local Lagrangian VE basis $\{\vvarphi_i\}_{i=1}^{\Ndof[E]}$  related to the local degrees of freedom, where the DOF numbering first counts the boundary DOFs and then the internal DOFs. Furthermore, for each element $E \in \Th$, we define the operators $\dof_i : \Vh[E]{k} \to \R$ which associate each function $\vv \in \Vh[E]{k}$ to its $i$-th degree of freedom.

Now, let us introduce the $L^2(E)$-projector $\vproj{0,E}{k}: \Vh[E]{k} \to \VPoly{k}{E}$, which is defined by the orthogonality condition
\begin{equation}
\scal[E]{\vv - \vproj{0,E}{k} \vv}{\qq} = 0 \quad \forall \qq \in \VPoly{k}{E},\ \vv \in \Vh[E]{k}.
\label{eq:proj}
\end{equation}
We note that each combination of the aforementioned degrees of freedom makes the projection $\vproj{0,E}{k} \vv_h$ of a function $\vv_h \in \Vh[E]{k}$ computable. In particular, the computation of $\vproj{0,E}{k} \vv_h$ with the pairs \ref{DOF:1a}-\ref{DOF:2a} and \ref{DOF:1a}-\ref{DOF:2b} has been largely discussed in \cite{MixedImplem, berrone2023orthogonal}. Concerning the choice \ref{DOF:1b}-\ref{DOF:2a}, we first note that, given $\vv \in \Vh[E]{k}$, the orthogonality condition \eqref{eq:proj} yields

\begin{multline}
 \scal[E]{\vproj{0,E}{k}  \vv}{\g^{\nabla,m}_{\alpha}} = \scal[E]{\vv}{\g^{\nabla,m}_{\alpha}} = \int_E \vv \cdot \nabla m_{\alpha + 1}^{k+1}\\
= - \int_E \div \vv \ m_{\alpha + 1}^{k+1} + \sum_{e \in \Eh[E]} \int_e \vv \cdot \nn_e \gamma_e\left(m_{\alpha + 1}^{k+1}\right), \quad \forall \alpha = 1,\dots,n_{k}^{\nabla} 
\label{eq:proj_nabla}
\end{multline}
and
\begin{equation}
\scal[E]{\vproj{0,E}{k} \vv}{\g^{\perp,m}_{\alpha}} = \scal[E]{\vv}{\g^{\perp,m}_{\alpha}}, \quad \forall \alpha = 1,\dots,n_{k}^{\perp},
\label{eq:proj_perp}
\end{equation}
where $\gamma_e\left(m_{\alpha + 1}^{k+1}\right)$ is the trace of the monomial $m_{\alpha + 1}^{k+1}$ of order $k+1$ on the edge $e \in \Eh[E]$.
Equation \eqref{eq:proj_perp} can be computed throughout the internal degrees of freedom \eqref{eq:idof_perp_mon}. Now, we recall that $\div \vv$ is a polynomial $\sum_{\alpha=1}^{n_k} c_{\alpha} m_{\alpha}^k \in \Poly{k}{E}$ whose coefficients $\{c_{\alpha}\}_{\alpha=1}^{n_k}$ can be determined by imposing
\begin{equation}
\int_E \div \vv\  m_{\beta}^k = \sum_{\alpha=1}^{n_k} c_{\alpha} \int_{E} m_{\alpha}^k m_{\beta}^k = - \int_E \vv \cdot \nabla m_{\beta}^k + \sum_{e \in \Eh[E]} \int_e \vv \cdot \nn_e \gamma_e\left(m_{\beta}^k\right), \quad \forall \beta = 1,\dots,n_k.
\label{eq:diver}
\end{equation}
The first term of the right-hand side of \eqref{eq:diver} can be computed throughout the internal degrees of freedom \eqref{eq:idof_nabla_mon}. Furthermore, we can write the trace of monomials as
\begin{equation}
\gamma_e\left(m_{\beta}^{k}\right) = \sum_{j=1}^{k+1} \mathbf{C}^e_{\beta j} t_j
\end{equation}
and compute the second term of the right-hand side of \eqref{eq:diver} 
by resorting to the boundary degrees of freedom \ref{DOF:1b} simply as
\begin{equation}
\int_e \vv \cdot \nn_e \gamma_e\left(m_{\beta}^k\right) = \sum_{j=1}^{k+1} \mathbf{C}^e_{\beta j} \int_0^{1} \widehat{\vv \cdot \nn_e} t_j \vert e \vert.
\end{equation}
In order to compute the second term of the right-hand side of equation \eqref{eq:proj_nabla}, we should determine the polynomial $\widehat{\vv \cdot \nn_e}$ on each edge $e \in \Eh[E]$.
However, if $\{\vvarphi_i^e\}_{i=1}^{k+1}$ is the local Lagrangian mixed VE basis related to the boundary degrees of freedom defined on the edge $e\in \Eh[E]$, we observe that
\begin{equation}
\widehat{\vvarphi_i^e \cdot \nn_e} = \frac{t_i}{\vert e \vert},\quad \forall i =1,\dots,k+1,
\end{equation}
while $\vvarphi \cdot \nn_e$ is the zero-polynomial if it is related to an internal degree of freedom or to a different edge of $E$.
Finally, since $\QPoly{k+1}{[0,1]}$ is an $L^2([0,1])$-orthonormal basis for $\Poly{k+1}{[0,1]}$, we simply have $\forall \alpha =1,\dots,n_k^{\nabla}$, $i = 1,\dots,k+1$ and $\forall e \in \Eh[E]$
\begin{equation}
\int_e \vvarphi_i^e \cdot \nn_e \gamma_e\left(m_{\alpha 
+ 1}^{k+1}\right) = \sum_{j=1}^{k+2} \mathbf{C}^e_{\beta j} \int_0^{1} \widehat{\vvarphi_i^e \cdot \nn_e} t_j \vert e \vert = \sum_{j=1}^{k+2} \mathbf{C}^e_{\alpha +1, j} \int_0^{1} t_i t_j = \mathbf{C}^e_{\alpha +1, i} \delta_{ij}.
\end{equation}

The construction of the method with the choice \ref{DOF:1b}-\ref{DOF:2b} is analogous to the one which exploits the degrees of freedom \ref{DOF:1b}-\ref{DOF:2a}. Indeed, we recall that we are able to write
\begin{linenomath}
\begin{equation*}
\g^{\nabla,\overline{q}}_{\alpha} = \sum_{\beta=1}^{n_k^{\nabla}} \mathbf{L}^{\nabla,k}_{\alpha \beta} \nabla q_{\beta+1} = \sum_{\beta=1}^{n_k^{\nabla}} \sum_{\gamma=1}^{n_{k+1}}\mathbf{L}^{\nabla}_{\alpha \beta} \mathbf{L}^{k+1}_{\beta+1, \gamma} \nabla m^{k+1}_{\gamma},
\end{equation*}
\end{linenomath}
where $\mathbf{L}^{\nabla,k}$ and $\mathbf{L}^{k+1}$ are defined in \eqref{eq:gnablaortho} and \eqref{eq:qbasis}, respectively.

\begin{remark}
Note that, since we define the one-dimensional polynomial basis $\QPoly{k}{[0,1]}$ on the interval $[0,1]$, we must perform the orthogonalization process just once. Thus, the additional cost in taking an $L^2([0,1])$-orthonormal basis instead of the one-dimensional monomial basis is negligible and independent of the number of edges of the tessellation $\Th$.

\end{remark}

\subsection{The Mixed Virtual Element Formulation of the model problem}

On each element $E \in \Th$, let us define the continuous local bilinear form  

\begin{equation*}
\dbilin[E]{\uu}{\vv} = \scal[E]{\D^{-1}\uu}{\vv},\quad\forall \uu,\vv \in \V
\end{equation*}
and its discrete counterpart 
\begin{equation}
\dbilinh[E]{\uu_h}{\vv_h} = \dbilinhC[E]{\uu_h}{\vv_h} + \stab[E]{\left(\bm{I}-\vproj{0,E}{k}\right)\uu_h}{\left(\bm{I}-\vproj{0,E}{k}\right)\vv_h}
\end{equation}
which is the sum of the consistency term

\begin{equation*}
 \dbilinhC[E]{\uu_h}{\vv_h} = \scal[E]{\D^{-1}\vproj{0,E}{k} \uu_h}{\vproj{0,E}{k}\vv_h}  
\end{equation*}
and of the \textit{stability term} $\stab[E]{}{}$, which is any symmetric positive definite bilinear form that satisfies
\begin{equation}
\alpha_{\ast} \dbilin[E]{\vv}{\vv} \leq \stab[E]{\vv}{\vv} \leq \alpha^{\ast} \dbilin[E]{\vv}{\vv}, \quad \forall \vv \in \Vh[E]{k}
\label{eq:stabprop}
\end{equation}
for some positive constants $\alpha_{\ast},\ \alpha^{\ast}$ depending on $\D^{-1}$ but independent on $h$ \cite{basicMixed, Mixed}.

Now, let us introduce the global mixed virtual element spaces

\begin{equation*}
\Vh{k} = \Big\{ \vv \in H_{0,\Gamma_N}(\div; \Omega): \ \vv_{|E} \in \Vh[E]{k}\ \forall E \in \Th\Big\},
\end{equation*}
\begin{equation*}
\Qh{k} = \Big\{ q \in L^2(\Omega): \ q_{|E} \in \Qh[E]{k} = \Poly{k}{E}\ \forall E \in \Th\Big\}.
\end{equation*}
for the velocity and the pressure variables, respectively. In particular, as global degrees of freedom for each $\vv_h \in \Vh{k}$, we consider
\begin{itemize}
\item the boundary degrees of freedom of $\vv_h$ defined on each internal edge of the decomposition and at edge boundary with Dirichlet boundary conditions;
\item the internal degrees of freedom in each element $E\in \Th$.
\end{itemize}
Furthermore, the value of the boundary DOFs at the Neumann edge is fixed in accordance with the value of the Neumann boundary conditions.

Finally, the virtual element discretization of the problem \eqref{eq:mixedForm} reads
\begin{equation}
\begin{cases}
\text{Find } (\uu_{0,h},p_h)\in \Vh{k} \times \Qh{k} \text{ such that } \uu_{h} = \uu_{0,h} + \uu_{N,h} \text{ and } p_h \text{ satisfy}  \\
\displaystyle\sum_{E \in \Th} (\dbilinh[E]{\uu_h}{\vv_h} - \scal[E]{p_h}{\div \vv_h}) = - \displaystyle\sum_{ E\in\Th} \displaystyle\sum_{\substack{e \in \Eh[E]:\\ e \subset \Gamma_D}}\langle g_D, \vv_h \cdot \nn_{e} \rangle_{\pm\frac{1}{2},e} &\forall \vv_h \in \Vh{k}\\
\displaystyle\sum_{E \in \Th} \scal[E]{\div \uu_h}{q_h} = \sum_{E \in \Th} \scal[E]{f}{q_h} & \forall q_h\in \Qh{k}
\end{cases},
\label{eq:mixedvemproblem}
\end{equation}
where $\uu_{N,h} \in \Big\{\vv \in H(\div;\Omega): \vv \in \Vh[E]{k} \forall E \in \Th \Big\}$ is such that $\dof_i(\uu_{N,h}) = \dof_i(\uu_N)$ for each boundary degree of freedom $i$.

\section{The stabilization term}\label{sec:stabilization}

Let us introduce the elemental matrix $\A^{E} \in \R^{\Ndof[E] \times \Ndof[E]}$, whose entries are defined as the application of the local discrete bilinear form $\dbilinh[E]{}{}$ to the Lagrangian basis functions of $\Vh[E]{k}$, i.e. $\forall i,j = 1,\dots, \Ndof[E]$

\begin{align*}
\left(\A^{E}\right)_{ij} &=  \dbilinh[E]{\vvarphi_i}{\vvarphi_j} \\
&= \dbilinhC[E]{\vvarphi_i}{\vvarphi_j} +  \stab[E]{(I - \vproj{0,E}{k})\vvarphi_i}{(I - \vproj{0,E}{k})\vvarphi_j}\\
& \coloneqq  \left(\A^{E}_C\right)_{ij} + \left(\A^{E}_S\right)_{ij},
\end{align*}
where $\A^{E}_C$ and $\A^{E}_S$ represent the elemental matrices related to the consistency and the stability term, respectively.
The complete elemental matrix related to the mixed discretization of the problem \eqref{eq:mixedvemproblem} reads 

\begin{equation*}
\KK^{E} = \begin{bmatrix}
\A^{E} & -(\mathbf{W}^E)^T\\
\mathbf{W}^E & \mathbf{0}
\end{bmatrix} \in \R^{(\Ndof[E]+n_k) \times (\Ndof[E]+n_k)},
\end{equation*}
where the entries of the divergence matrix $\mathbf{W}^E \in \R^{n_k \times \Ndof[E]}$ are defined as

\begin{equation*}
\mathbf{W}^E _{\alpha i } = \scal[E]{p_{\alpha}}{\vvarphi_i},\quad \forall p_{\alpha} \in \M{k}{E} ( \text{or } \forall p_{\alpha} \in \QPoly{k}{E}),\ \forall i =1,\dots,\Ndof[E]. 
\end{equation*}
Since the degrees of freedom of the velocity space are chosen in such a way the related Lagrangian VE basis functions scale uniformly with respect to the mesh size $h$, the most natural mixed VEM stabilization $\stab[E]{}{}$ which satisfies \eqref{eq:stabprop} is the so-called \textit{dofi-dofi} stabilization \cite{LBe13, Mixed}:
\begin{equation}
\stabdof[E]{\uu- \vproj{0,E}{k} \uu}{\vv- \vproj{0,E}{k} \vv} = C_{\D^{-1}}\vert E \vert \sum_{i=1}^{\Ndof[E]} \dof_i(\uu- \vproj{0,E}{k} \uu) \dof_i(\vv- \vproj{0,E}{k} \vv),
\end{equation}
where $C_{\D^{-1}}$ is a constant depending on $\D^{-1}$. Moreover, since both the spaces $\GPoly{k-1}{\nabla}{E}$ and $\GPoly{k}{\perp}{E}$ represent polynomials in $\VPoly{k}{E}$, it follows
\begin{equation}
\dof_i(\uu- \vproj{0,E}{k} \uu) = 0
\label{eq:idof_stab}
\end{equation}
for each internal degree of freedom $i$. Thus, in the mixed VEM construction, it is not necessary to stabilize the internal degrees of freedom.

Furthermore, as highlighted in \cite{BEIRAODAVEIGA20171110}, in order to avoid to level off the stabilization term with respect to the consistency term for the higher polynomial degrees, which would lead to a loss of accuracy, we can choose the so-called \textit{D-recipe} stabilization, defined as follows
\begin{equation}
\stabD[E]{\uu- \vproj{0,E}{k} \uu}{\vv- \vproj{0,E}{k} \vv} = \sum_{i=1}^{\Ndof[E]} S_{ii} \dof_i(\uu- \vproj{0,E}{k} \uu) \dof_i(\vv- \vproj{0,E}{k} \vv),
\label{eq:stab_drecipe}
\end{equation}
where $S_{ii} = C_{\D^{-1}} \vert E \vert \max(1, (\A^E_C)_{ii})$ if $i$ is related to a boundary degree of freedom and $S_{ii} = 0$ otherwise, since we do not need to stabilize the internal degrees of freedom (see equation \eqref{eq:idof_stab}). 

Usually, the constant $C_{\D^{-1}}$ is taken equal to the spectral norm $\| \D^{-1}\| = 1/\D_{min}$, since $\D$ is assumed to be symmetric and strong elliptic. 

Finally, the choice of the stabilization term and, in particular, of the constant $C_{\D^{-1}}$ should be dependent on the problem features and on the definition of the local degrees of freedom \cite{LBe13, russo2023quantitative}.

\section{Numerical experiments}\label{sec:experiments}

In this section, we perform some numerical experiments that allow us to show the role of the boundary degrees of freedom and of the stabilization term in preventing the ill-conditioning of the system matrix. 
To this end, we analyze the behaviour of the global system matrix $\KK$ and of the following errors:
\begin{equation}
\mathrm{err}_p = \frac{\sqrt{\sum_{E\in\Th}\norm[E]{p-p_h}^2}}{\norm[\Omega]{p}}
\label{eq:errorp}
\end{equation}
\begin{equation}
\mathrm{err}_{\uu} = \frac{\sqrt{\sum_{E\in\Th}\norm[E]{\uu-\vproj{0,E}{k}\uu_h}^2}}{\norm[\Omega]{\uu}}
\label{eq:errorv}
\end{equation}
at varying of the polynomial degree $k$ or of the mesh size $h$, for different families of meshes.
Given $k \geq 0$ and the mesh size $h$, we recall that if the solution is sufficiently smooth, the expected convergence rates of errors \eqref{eq:errorp} and \eqref{eq:errorv} is $O(h^{k+1})$.

In the following, we use the notation
\begin{itemize}
\item Mon (a) to denote the approach which exploits the pair of DOFs \ref{DOF:1a}-\ref{DOF:2a};
\item Mon (b) to denote the approach which exploits the pair of DOFs \ref{DOF:1b}-\ref{DOF:2a};
\item Ortho (a) to denote the approach which exploits the pair of DOFs \ref{DOF:1a}-\ref{DOF:2b};
\item Ortho (b) to denote the approach which exploits the pair of DOFs \ref{DOF:1b}-\ref{DOF:2b}.
\end{itemize}
We note that in the monomial approaches (Mon), we use the scaled monomial basis as the basis for the pressure space, while in the orthonormal approaches (Ortho), we use the $\QPoly{k}{E}$ basis as the polynomial basis.

\subsection{Test 1: Boundary degrees of freedom}

In this first test, we analyze the behaviour of the four aforementioned approaches by solving a Poisson problem with homogeneous Dirichlet boundary conditions.

More precisely, let us set $\Omega = (0,2)^2$ and we define the forcing term $f$ in such a way the exact pressure is

\begin{equation*}
p(x,y) = \sin(\pi x) \sin(\pi y).
\end{equation*} 

\begin{figure}[]
	\centering
	\subfigure[\label{fig:M1}]
	{\includegraphics[width=.32\textwidth, height = .22\textheight]{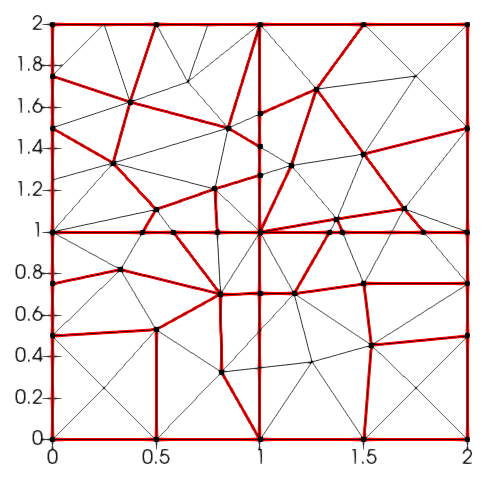}}
	\subfigure[\label{fig:M2}]
	{\includegraphics[width=.32\textwidth, height = .22\textheight]{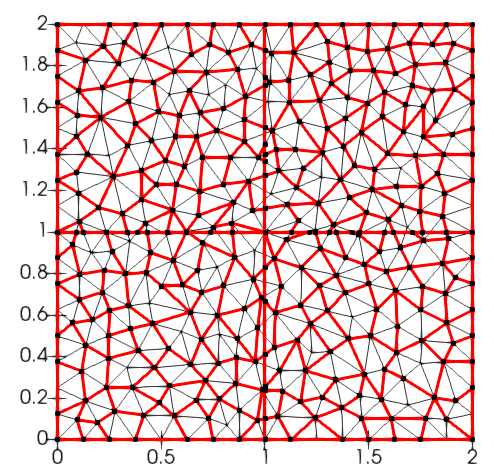}}
	\subfigure[\label{fig:M3}]
	{\includegraphics[width=.32\textwidth, height = .215\textheight]{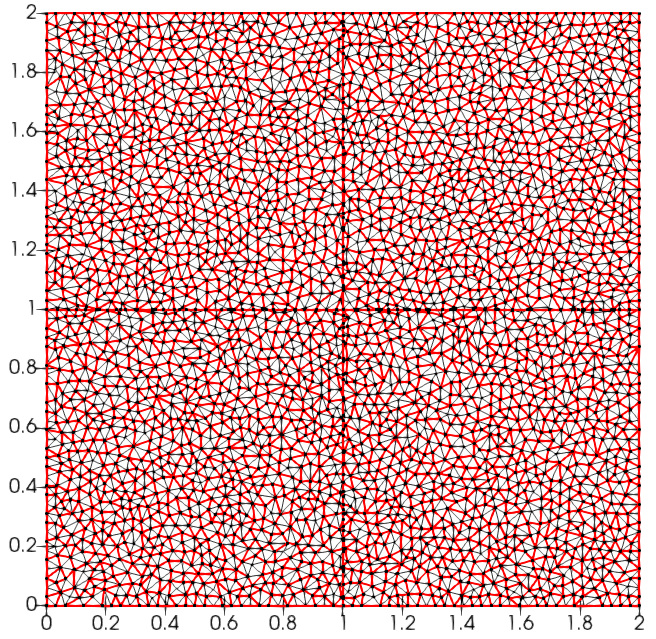}}
	\caption{Test 1. The three concave refinements $\Th[i]^C$, $i=1,2,3$.}
	\label{fig:mesh_test1}
\end{figure}

In this test, we employ the dofi-dofi stabilization term with $C_{\D^{-1}} = 1$ and we evaluate the performances of our approaches on a family of three concave meshes $\{\Th[i]^C\}_{i=1}^3$ which are generated throughout an agglomeration process starting from triangular meshes with a different refinement level, as shown in Figure \ref{fig:mesh_test1}. 

\begin{figure}[]
	\centering
	\subfigure[\label{fig:condStiff_M1_Poisson}]
	{\includegraphics[width=.32\textwidth, height = .25\textheight]{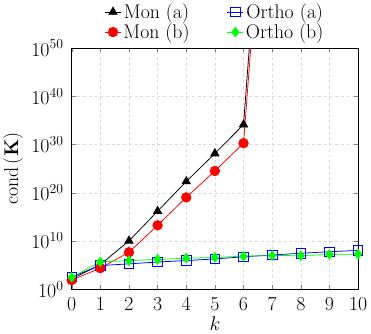}}
	\subfigure[\label{fig:condStiff_M2_Poisson}]
	{\includegraphics[width=.32\textwidth, height = .25\textheight]{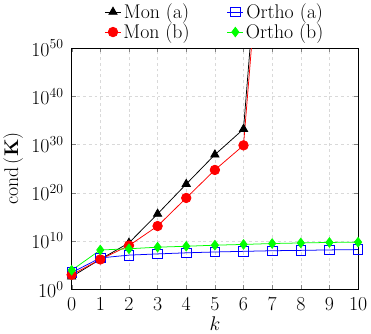}}
	\subfigure[\label{fig:condStiff_M3_Poisson}]
	{\includegraphics[width=.32\textwidth, height = .25\textheight]{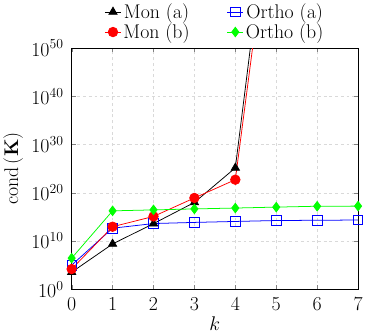}}
	\caption{Test 1. Condition number of $\mathbf{K}$ vs. $k$. Left: Mesh $\Th[1]^C$. Center: Mesh $\Th[2]^C$. Right: Mesh $\Th[3]^C$.}
	\label{fig:condStiff_test1}

	\centering
	\subfigure[\label{fig:ErrorL2Pressure_M1_Poisson}]
	{\includegraphics[width=.32\textwidth, height = .25\textheight]{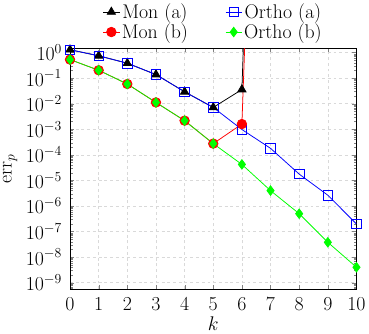}}
	\subfigure[\label{fig:ErrorL2Pressure_M2_Poisson}]
	{\includegraphics[width=.32\textwidth, height = .25\textheight]{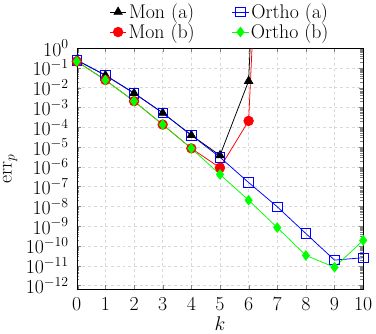}}
	\subfigure[\label{fig:ErrorL2Pressure_M3_Poisson}]
	{\includegraphics[width=.32\textwidth, height = .25\textheight]{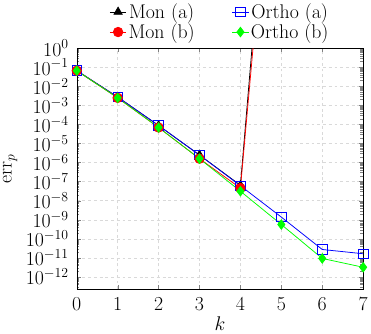}}
	\caption{Test 1. Behaviour of $\mathrm{err}_p$ \eqref{eq:errorp} vs. $k$. Left: Mesh $\Th[1]^C$. Center: Mesh $\Th[2]^C$. Right: Mesh $\Th[3]^C$.}
	\label{fig:errorL2Pressure_test1}
	
	\centering
	\subfigure[\label{fig:ErrorL2Velocity_M1_Poisson}]
	{\includegraphics[width=.32\textwidth, height = .25\textheight]{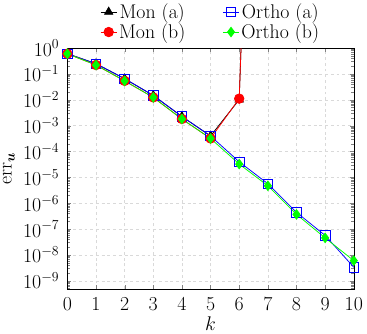}}
	\subfigure[\label{fig:ErrorL2Velocity_M2_Poisson}]
	{\includegraphics[width=.32\textwidth, height = .25\textheight]{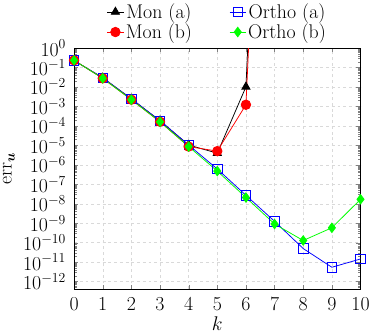}}
	\subfigure[\label{fig:ErrorL2Velocity_M3_Poisson}]
	{\includegraphics[width=.32\textwidth, height = .25\textheight]{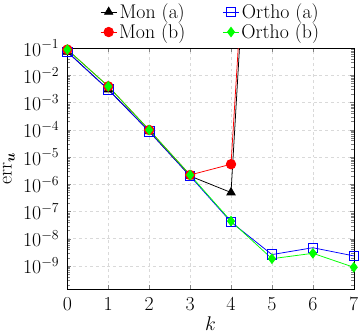}}

	\caption{Test 1. Behaviour of $\mathrm{err}_{\uu}$ \eqref{eq:errorv} vs. $k$. Left: Mesh $\Th[1]^C$. Center: Mesh $\Th[2]^C$. Right: Mesh $\Th[3]^C$.}
	\label{fig:errorL2Velocity_test1}
\end{figure}

In Figure \ref{fig:condStiff_test1}, we show the behaviour of the condition number of the global system matrix $\KK$ at varying of the polynomial degree $k$, for each concave mesh $\Th[i]^C$, $i=1,2,3$, in semilog plots. From these graphs, we note that changing the boundary degrees of freedom from \ref{DOF:1a} to \ref{DOF:1b} generally does not ensure an improvement in the condition number of the global system matrix for fixed internal degrees of freedom. Furthermore, we observe that, in order to cure the ill-conditioning of the global system matrix, the use of an $L^2(E)$-orthonormal (vector) polynomial basis for $\VPoly{k}{E}$ is strongly recommended, as already highlighted in \cite{berrone2023orthogonal}.

\begin{figure}[]
	\centering
	\subfigure[\label{fig:ErrorHPressure_1}]
	{\includegraphics[width=.32\textwidth, height = .25\textheight]{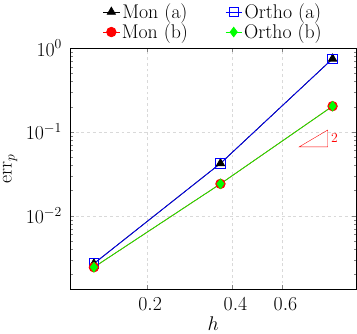}}
	\subfigure[\label{fig:ErrorHPressure_3}]
	{\includegraphics[width=.32\textwidth, height = .25\textheight]{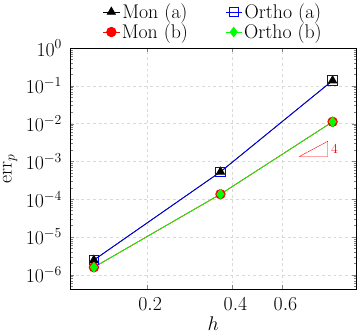}}
	\subfigure[\label{fig:ErrorHPressure_5}]
	{\includegraphics[width=.32\textwidth, height = .25\textheight]{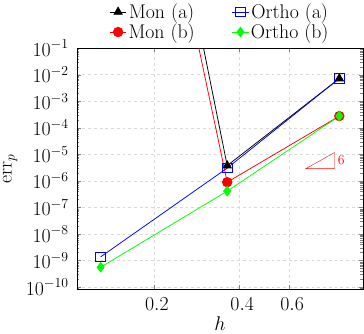}}
	\caption{Test 1. Behaviour of $\mathrm{err}_p$ \eqref{eq:errorp} vs. $h$. Left: $k=1$. Center: $k=3$. Right: $k=5$.}
	\label{fig:errorH_L2Pressure_test1}

	\centering
	\subfigure[\label{fig:ErrorHVelocity_1}]
	{\includegraphics[width=.32\textwidth, height = .25\textheight]{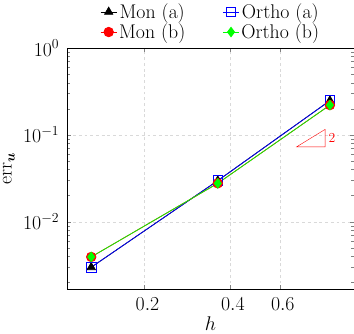}}
	\subfigure[\label{fig:ErrorHVelocity_3}]
	{\includegraphics[width=.32\textwidth, height = .25\textheight]{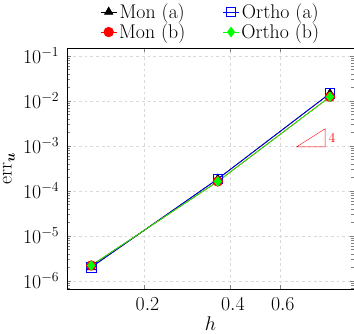}}
	\subfigure[\label{fig:ErrorHVelocity_5}]
	{\includegraphics[width=.32\textwidth, height = .25\textheight]{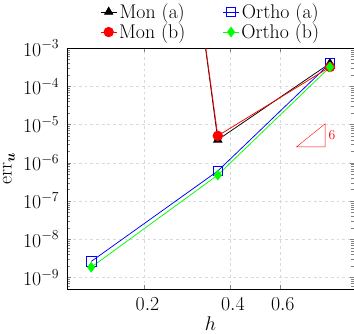}}
	\caption{Test 1. Behaviour of $\mathrm{err}_{\uu}$ \eqref{eq:errorv} vs. $h$. Left: $k=1$. Center: $k=3$. Right: $k=5$.}
	\label{fig:errorH_L2L2Velocity_test1}
\end{figure}

Figures \ref{fig:errorL2Pressure_test1} and \ref{fig:errorL2Velocity_test1} show the behaviour of errors \eqref{eq:errorp} and \eqref{eq:errorv} at varying of the polynomial degree $k$ for each $\Th[i]^C$, with $i=1,2,3$, in semilog plots. Furthermore, Figures \ref{fig:errorH_L2Pressure_test1} and \ref{fig:errorH_L2L2Velocity_test1} show the behaviour of such errors for decreasing values of the mesh size $h_i$, $i=1,2,3$, for $k=1,3,5$, with a loglog scale. From these figures, we can note that changing the internal degrees of freedom from \ref{DOF:2a} to \ref{DOF:2b} does not modify significantly the behaviour of errors \eqref{eq:errorp} and \eqref{eq:errorv}, at varying of the mesh size $h$, for the lower values of the polynomial degree $k$. In general, this is not true for the boundary degrees of freedom. Indeed, from Figures \ref{fig:errorL2Pressure_test1} and \ref{fig:errorH_L2Pressure_test1}, we can note that the error \eqref{eq:errorp} is sensitive to a variation from \ref{DOF:1a} to \ref{DOF:1b} of the boundary degrees of freedom, especially on the coarser meshes. As the mesh is refined, such difference becomes smaller and smaller and the orthonormal approaches tend to behave in the same way regardless of the type of boundary DOFs used. 

Finally, for the higher values of $k$, the errors start to raise due to the ill-conditioning of the matrix $\KK$ in the Mon approaches, while the Ortho approaches are robust also for the higher polynomial degrees. 

\subsection{Test 2: Anisotropic diffusion problems}

In this experiment, we want to analyze the sensitivity of the presented approaches to the choice of stabilization in a context where such sensitivity becomes the main issue to overcome, namely the diffusion problems with high anisotropic coefficients.

Equations characterized by anisotropic diffusion coefficients arise in many practical contexts, such as the heat equation, groundwater flow, transport problems and so on. Generally, these types of problems are expressed as parametric problems and they are numerically treated by means of ad hoc methods, needed to avoid the so-called \textit{locking} phenomenon \cite{locking}. This phenomenon occurs experimentally when the discretization error does not decrease at the expected rate when the parameter tends to limiting values and, in general, is typical of the lower order schemes.
These ad hoc methods include \textit{variational crimes}, i.e. modification of the bilinear form \cite{Havu2021}, and flow-aligned grid methods \cite{alignment}.
In particular, in the Virtual Element context, the isotropic nature of the standard stabilization term can become an issue in these kinds of problems and different approaches have been studied to handle the anisotropic nature of the diffusion tensors \cite{berrone2023lowest, MAZZIA202063} mainly for the primal formulation of the method.

Thus, we consider the test problem proposed in \cite{MANZINI2007751}, which is a dimensionless parametric version of problem \eqref{eq:primalForm} with a constant diffusion tensor, defined on $\Omega = (0,1)^2$. In particular, the diffusion tensor $\D = \begin{bmatrix} 1 & 0 \\ 0 & \epsilon \end{bmatrix}$ depends on the diffusion parameter $\epsilon \in [10^{-6}, 1]$, which, in this case, represents also the anisotropic ratio. In our notation, $\D_{min} = \epsilon$ (or $\D^{-1}_{max} = \frac{1}{\epsilon}$) and $\D_{max} = 1$.

\begin{figure}[]
	\centering
	\subfigure[\label{fig:Q5x5}]{\includegraphics[width=.24\textwidth, height = .17\textheight]{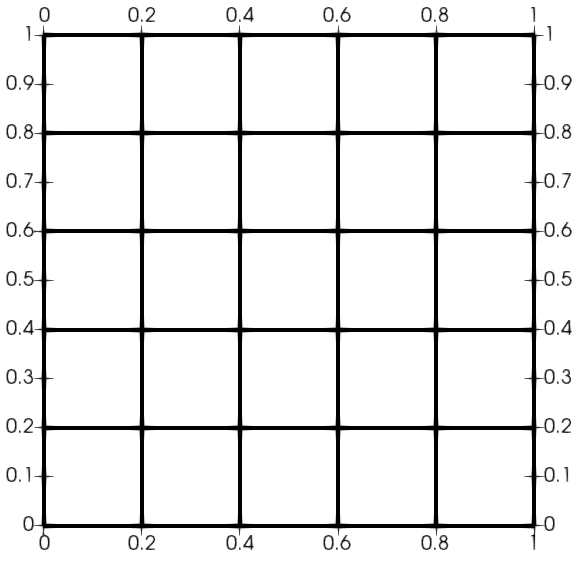}}
	\subfigure[\label{fig:QD5x5}]{\includegraphics[width=.24\textwidth, height = .17\textheight]{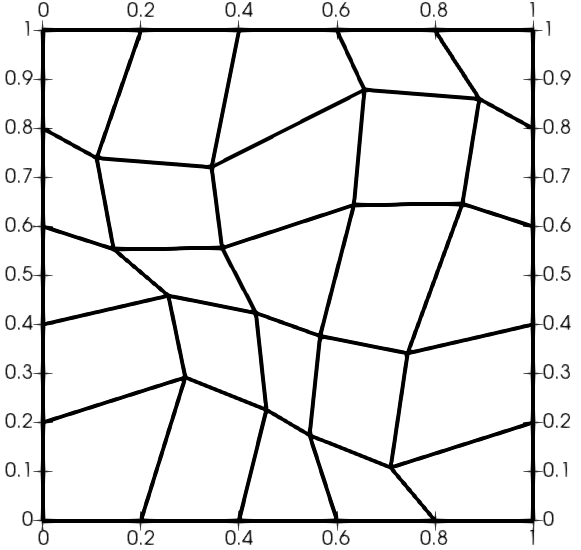}}
	\subfigure[\label{fig:Q40x40}]{\includegraphics[width=.24\textwidth, height = .17\textheight]{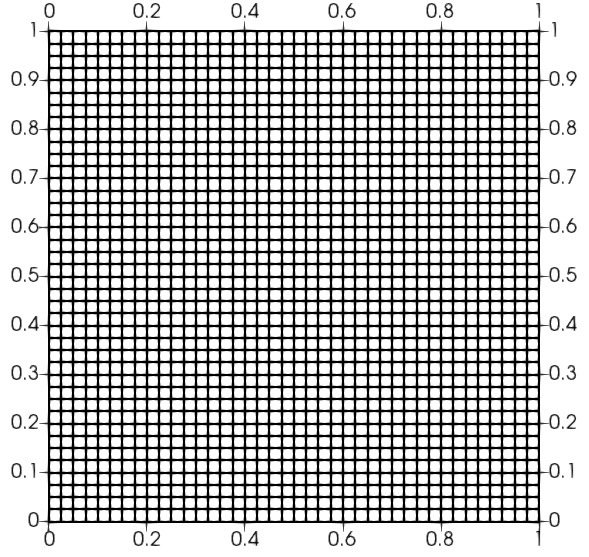}}
	\subfigure[\label{fig:QD40x40}]{\includegraphics[width=.24\textwidth, height = .17\textheight]{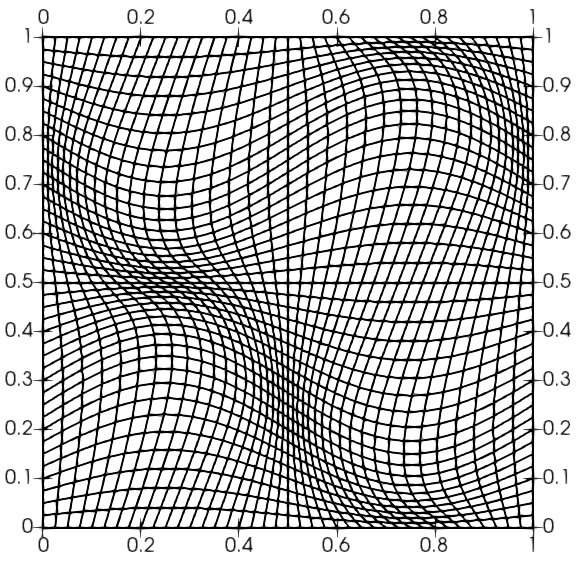}}
	\caption{Test 2. \ref{fig:Q5x5}: Mesh $\Th[1]^{Q}$. \ref{fig:QD5x5}: Mesh $\Th[1]^{DQ}$. \ref{fig:Q40x40}: Mesh $\Th[4]^{Q}$. \ref{fig:QD40x40}: Mesh $\Th[4]^{DQ}$.}
	\label{fig:mesh_test2}
\end{figure}

The performances of the four approaches are evaluated on two different kinds of families of meshes: a cartesian $\Th^Q$ family and a family $\Th^{DQ}$ of distorted quadrilateral meshes obtained by the cartesian ones throughout a sine distortion. For each family of meshes, we consider four refinements $\{\Th[i]^Q\}_{i=1}^{4}$ and $\{\Th[i]^{DQ}\}_{i=1}^{4}$: the first and the last refinement of each family are shown in Figure \ref{fig:mesh_test2}.

In  order to compute errors \eqref{eq:errorp} and \eqref{eq:errorv}, we choose the parametric exact solution
\begin{equation}
p(x,y) = \exp(-2 \pi \sqrt{\epsilon} x)\sin( 2 \pi y).
\label{eq:solution_test2}
\end{equation}
The presence of $\epsilon$ at the exponent of \eqref{eq:solution_test2} makes the low conductivity direction dominant when $\epsilon$ tends to zero and the nearly pure Neumann boundary conditions are set, by leading, in general, to very poor results when employing standard methods \cite{Havu2021}. Thus, we test three different kinds of boundary conditions (BCs in short): 
\begin{itemize}
\item pure Dirichlet boundary conditions, i.e. $\Gamma_D = \Gamma$;
\item mixed Dirichlet-Neumann boundary conditions with
\begin{linenomath}
\begin{equation*}
\Gamma_D = \{(x,y): x = 0  \text{ or } y = 0\};
\end{equation*}
\end{linenomath}
\item nearly pure Neumann conditions, that is we set
\begin{linenomath}
\begin{equation*}
\Gamma_D = \{(x,y): (x=1 \text{ and } 1-\delta \leq y \leq 1) \text{ or } (y=1 \text{ and } 1-\delta \leq x \leq 1)\}, 
\end{equation*}
\end{linenomath}
where $\delta$ decreases with the mesh size as $\frac{1}{5 \cdot 2^{i-1}}$ $i=1,\dots,4$.
\end{itemize}
In the first two cases, generally, no locking phenomenon occurs.

Furthermore, we test three possible choices for the stabilization term, namely
\begin{itemize}
\item S$_1$: the standard dofi-dofi stabilization with $C_{\D^{-1}} = \|\D^{-1} \| = \frac{1}{\epsilon}$;
\item S$_2$: the standard dofi-dofi stabilization with $C_{\D^{-1}} = 1$;
\item S$_3$: the D-recipe stabilization with $C_{\D^{-1}} = 1$.
\end{itemize}
We observe that when $\epsilon$ becomes very small, the constant $C_{\D^{-1}}$ related to the choice S$_1$ becomes very big.

\subsubsection{Effect of the anisotropy on the condition number of the global system matrix}

\begin{figure}[p]
\vspace{-0.5\baselineskip}
	\centering
\subfigure[\label{fig:condGlobStiff2D_Q5x5_BC1_Anisotropic_Ep_1.00000}]
	{\includegraphics[width=.35\textwidth, height = .20\textheight]{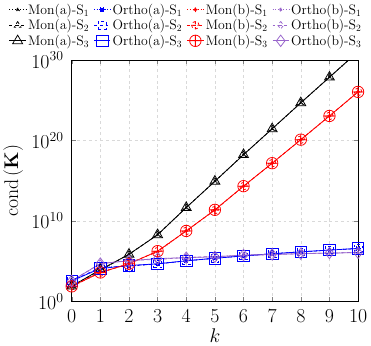}}
	\subfigure[\label{fig:condGlobStiff2D_Q5x5_BC1_Anisotropic_Ep_1em06}]
	{\includegraphics[width=.35\textwidth, height = .20\textheight]{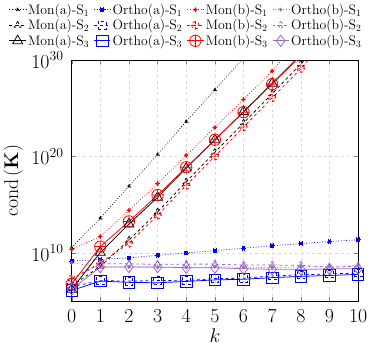}}
	
\vspace{-0.5\baselineskip}	\subfigure[\label{fig:condGlobStiff2D_QD5x5_BC1_Anisotropic_Ep_1}]
	{\includegraphics[width=.35\textwidth, height = .20\textheight]{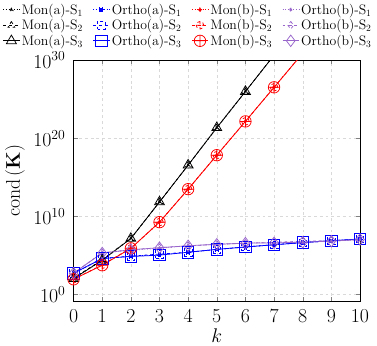}}
	{\includegraphics[width=.35\textwidth, height = .20\textheight]{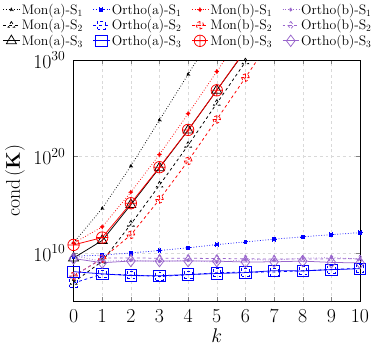}}
	\vspace{-0.5\baselineskip}
	\caption{Test 2. Condition number of $\KK$ vs. $k$. Left: $\epsilon = 1$. Right: $\epsilon = 10^{-6}$. 
First row: $\Th[1]^{Q}$. Second row: $\Th[1]^{DQ}$. Dirichlet BCs.}
	\label{fig:condStiff_test2_dirichlet}
	\centering
\vspace{-0.5\baselineskip}	\subfigure[\label{fig:condGlobStiff2D_Q5x5_BC3_Anisotropic_Ep_1}]
	{\includegraphics[width=.35\textwidth, height = .20\textheight]{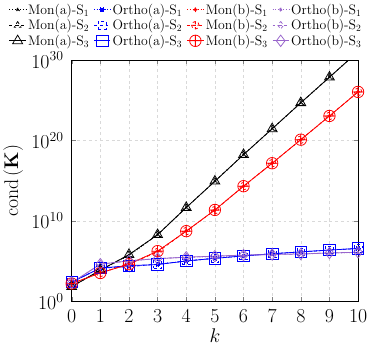}}
	\subfigure[\label{fig:condGlobStiff2D_Q5x5_BC3_Anisotropic_Ep_1em06}]
	{\includegraphics[width=.35\textwidth, height = .20\textheight]{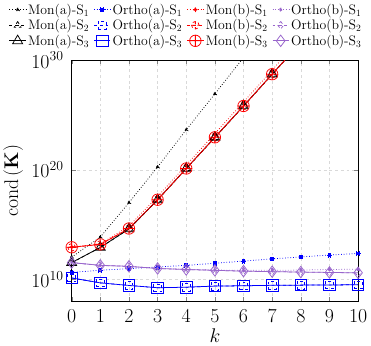}}

	\vspace{-0.5\baselineskip}	\subfigure[\label{fig:condGlobStiff2D_QD5x5_BC3_Anisotropic_Ep_1}]
	{\includegraphics[width=.35\textwidth, height = .20\textheight]{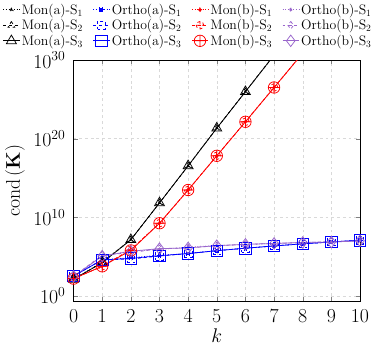}}
	\subfigure[\label{fig:condGlobStiff2D_QD5x5_BC3_Anisotropic_Ep_1em06}]
	{\includegraphics[width=.35\textwidth, height = .21\textheight]{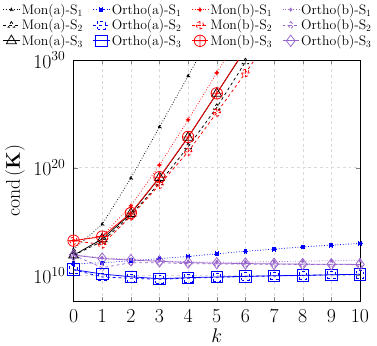}}
	\vspace{-0.5\baselineskip}
	\caption{Test 2. Condition number of $\KK$ vs. $k$. Left: $\epsilon = 1$. Right: $\epsilon = 10^{-6}$. 
First row: $\Th[1]^{Q}$. Second row: $\Th[1]^{DQ}$. Nearly pure Neumann BCs.}
	\label{fig:condStiff_test2_neumann}
\end{figure}

In Figures \ref{fig:condStiff_test2_dirichlet} and \ref{fig:condStiff_test2_neumann} we report the behaviour of the condition number of the global system matrix at varying of $k$ in semilog plots, when the Dirichlet and the nearly pure Neumann boundary conditions are set, respectively.
The results are related to $\epsilon \in \{1, 10^{-6}\}$ and to the $\Th[1]^{Q}$ and the $\Th[1]^{DQ}$ meshes.

Accordingly to results presented in \cite{berrone2023orthogonal}, we observe an exponential growth in the condition number of the matrix $\KK$ when the internal DOFs \ref{DOF:2a} are employed. A linear growth is observed instead when resorting to the choice \ref{DOF:2b}.  Furthermore, as already pointed out in the previous test, changing the boundary DOFs from \ref{DOF:1a} to \ref{DOF:1b} does not lead generally to an improvement of the behaviour of the condition number of $\KK$.

Furthermore, we note that a sine distortion of elements causes a faster increase in the condition number of $\KK$ when the internal DOFs \ref{DOF:2a} are used, while this growth is not so evident in the case of the internal DOFs \ref{DOF:2b}.

We further note that having nearly pure Neumann boundary conditions has just a small effect on the condition number of $\KK$ for the lower values of $k$ and that the condition number of $\KK$ seems to be mainly controlled by the anisotropic effect accordingly to what observed in \cite{MANZINI2007751}.
 
Finally, we observe that the pair \ref{DOF:1b}-\ref{DOF:2b} reveals to be the more robust approach with respect to the choice of the stabilization term, whereas stabilization choice S$_1$ seems to be the worst choice in terms of the condition number of $\KK$, if a combination of DOFs different from \ref{DOF:1b}-\ref{DOF:2b} is used.

\subsubsection{The mesh alignment and the locking phenomenon}

\begin{figure}[p]
	\centering
	\subfigure[\label{fig:ErrorL2Pressure_Q5x5_BC1_Anisotropic_Ep_1}]
	{\includegraphics[width=.35\textwidth, height = .26\textheight]{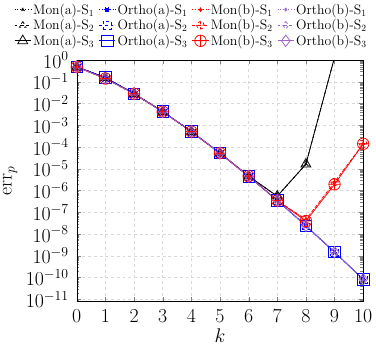}}
	\subfigure[\label{fig:ErrorL2Pressure_Q5x5_BC1_Anisotropic_Ep_1em06}]
	{\includegraphics[width=.35\textwidth, height = .25\textheight]{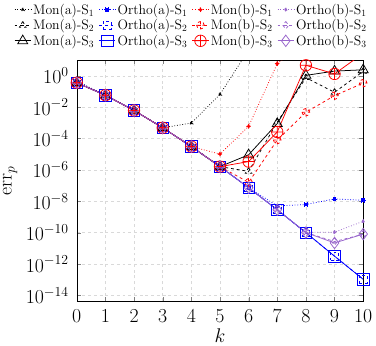}}
	\vspace{-0.5\baselineskip}
	\caption{Test 2. Behaviour of $\mathrm{err}_p$ \eqref{eq:errorp} vs. $k$, for $\Th[1]^{Q}$. Left: $\epsilon = 1$. Right: $\epsilon = 10^{-6}$. Dirichlet BCs.}
	\label{fig:errorL2Pressure_Q_test2_dirichlet}
	\centering
	\subfigure[\label{fig:ErrorL2Pressure_Q5x5_BC2_Anisotropic_Ep_1}]
	{\includegraphics[width=.35\textwidth, height = .25\textheight]{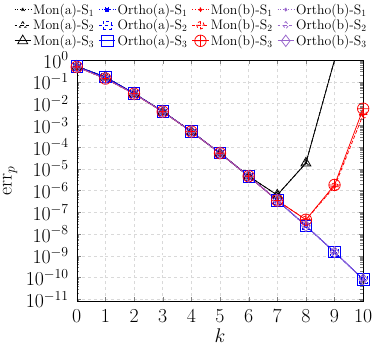}}
	\subfigure[\label{fig:ErrorL2Pressure_Q5x5_BC2_Anisotropic_Ep_1em06}]
	{\includegraphics[width=.35\textwidth, height = .25\textheight]{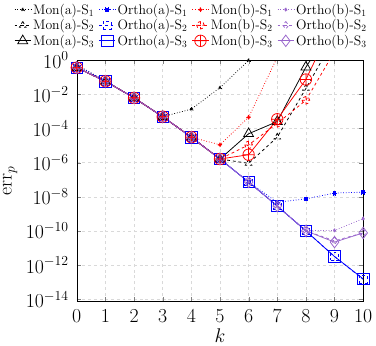}}
	\vspace{-0.5\baselineskip}
	\caption{Test 2. Behaviour of $\mathrm{err}_p$ \eqref{eq:errorp} vs. $k$, for $\Th[1]^{Q}$. Left: $\epsilon = 1$. Right: $\epsilon = 10^{-6}$. Mixed BCs.}
	\label{fig:errorL2Pressure_Q_test2_mixed}	
	\centering
	\subfigure[\label{fig:ErrorL2Pressure_Q5x5_BC3_Anisotropic_Ep_1}]
	{\includegraphics[width=.35\textwidth, height = .25\textheight]{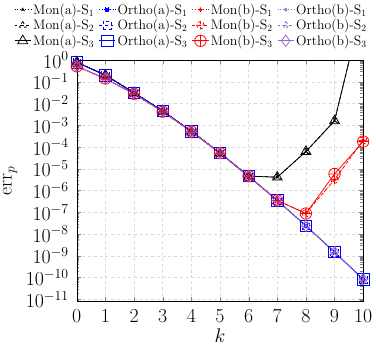}}
	\subfigure[\label{fig:ErrorL2Pressure_Q5x5_BC3_Anisotropic_Ep_1em06}]
	{\includegraphics[width=.35\textwidth, height = .25\textheight]{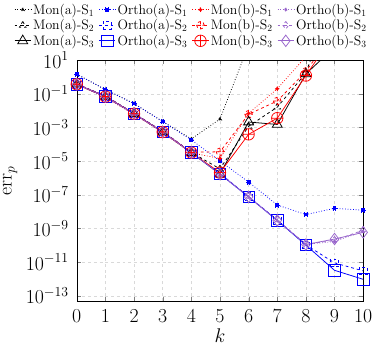}}
	\vspace{-0.5\baselineskip}	
	\caption{Test 2. Behaviour of $\mathrm{err}_p$ \eqref{eq:errorp} vs. $k$, for $\Th[1]^{Q}$. Left: $\epsilon = 1$. Right: $\epsilon = 10^{-6}$. Nearly pure Neumann BCs.}
	\label{fig:errorL2Pressure_Q_test2_neumann}
\end{figure}

Figures \ref{fig:errorL2Pressure_Q_test2_dirichlet}, \ref{fig:errorL2Pressure_Q_test2_mixed} and \ref{fig:errorL2Pressure_Q_test2_neumann} show the behaviour of the pressure error \eqref{eq:errorp} at varying of the polynomial degree $k$ related to the Dirichlet, mixed and nearly pure Neumann boundary conditions, respectively. These results are obtained on the cartesian mesh $\Th[1]^{Q}$.

From these figures we observe that, after an initial decrease, the error starts to raise due to ill-conditioning, but only when the internal DOFs \ref{DOF:2a} are employed. Choosing internal DOFs \ref{DOF:2b} leads to the best performances in each tested case for the higher values of the polynomial degree $k$.

The error curves related to the different analyzed approaches are very similar for the lower values of $k$ when a cartesian mesh is used. The only exception is represented by the choice boundary DOFs \ref{DOF:1a} and stabilization term S$_1$. In this case, error curves are slightly upward shifted for the smaller values of $\epsilon$ when nearly pure Neumann boundary conditions are set.

\begin{figure}[p]
	\centering	
		\subfigure[\label{fig:ErrorL2Pressure_QD5x5_BC1_Anisotropic_Ep_1}]
	{\includegraphics[width=.35\textwidth, height = .25\textheight]{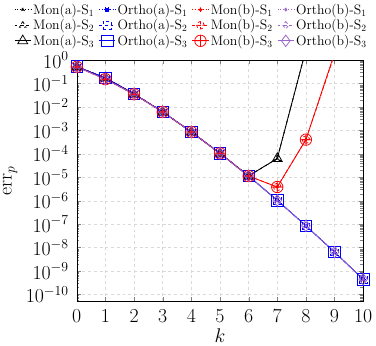}}
	\subfigure[\label{fig:ErrorL2Pressure_QD5x5_BC1_Anisotropic_Ep_1em06}]
	{\includegraphics[width=.35\textwidth, height = .25\textheight]{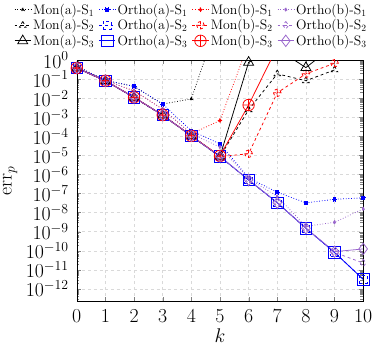}}
	\vspace{-0.5\baselineskip}
	\caption{Test 2. Behaviour of $\mathrm{err}_p$ \eqref{eq:errorp} vs. $k$, for $\Th[1]^{DQ}$. Left: $\epsilon = 1$. Right: $\epsilon = 10^{-6}$. Dirichlet BCs.}
	\label{fig:errorL2Pressure_DQ_test2_dirichlet}

	\centering	
	\subfigure[\label{fig:ErrorL2Pressure_QD5x5_BC2_Anisotropic_Ep_1}]
	{\includegraphics[width=.35\textwidth, height = .25\textheight]{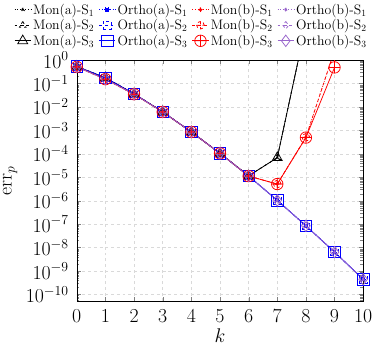}}
	\subfigure[\label{fig:ErrorL2Pressure_QD5x5_BC2_Anisotropic_Ep_1em06}]
	{\includegraphics[width=.35\textwidth, height = .25\textheight]{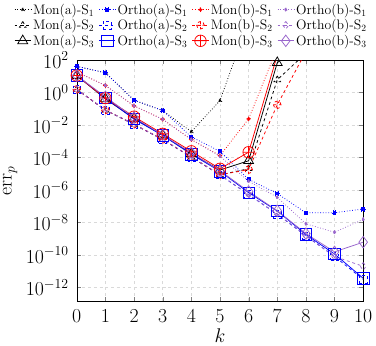}}
	\vspace{-0.5\baselineskip}
	\caption{Test 2. Behaviour of $\mathrm{err}_p$ \eqref{eq:errorp} vs. $k$, for $\Th[1]^{DQ}$. Left: $\epsilon = 1$. Right: $\epsilon = 10^{-6}$. Mixed BCs.}
	\label{fig:errorL2Pressure_DQ_test2_mixed}

	\centering	
		\subfigure[\label{fig:ErrorL2Pressure_QD5x5_BC3_Anisotropic_Ep_1}]
	{\includegraphics[width=.35\textwidth, height = .25\textheight]{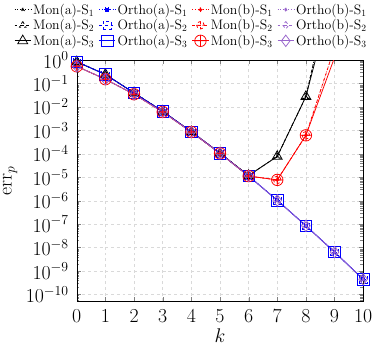}}
	\subfigure[\label{fig:ErrorL2Pressure_QD5x5_BC3_Anisotropic_Ep_1em06}]
	{\includegraphics[width=.35\textwidth, height = .25\textheight]{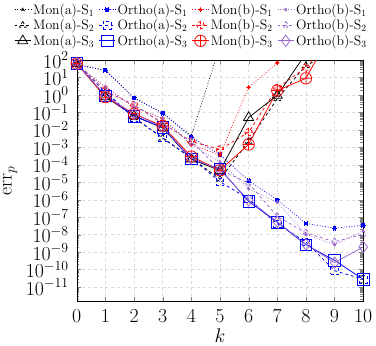}}
	\vspace{-0.5\baselineskip}
	\caption{Test 2. Behaviour of $\mathrm{err}_p$ \eqref{eq:errorp} vs. $k$, for $\Th[1]^{DQ}$. Left: $\epsilon = 1$. Right: $\epsilon = 10^{-6}$. Nearly pure Neumann BCs.}
	\label{fig:errorL2Pressure_DQ_test2_neumann}
\end{figure}

In Figures \ref{fig:errorL2Pressure_DQ_test2_dirichlet}, \ref{fig:errorL2Pressure_DQ_test2_mixed} and \ref{fig:errorL2Pressure_DQ_test2_neumann} we report the behaviour of the pressure error \eqref{eq:errorp} at varying of $k$ related to the Dirichlet, mixed and nearly pure Neumann boundary conditions, in the case of the distorted cartesian mesh $\Th[1]^{DQ}$.

By comparing these results with those obtained in the case of cartesian mesh, we can observe that, in the case of distorted meshes, the considered approaches show very different behaviours in terms of error \eqref{eq:errorp} when $\epsilon$ is very small. 
Indeed, we highlight that the cartesian mesh is aligned with the directions of the anisotropy, by limiting the effect of anisotropy. 
The main variations are observed for the approaches that exploit the stabilization term S$_1$ also in the case of distorted meshes.
Furthermore, we must observe an initial upward shift of the error curves related to the D-recipe S$_3$  for the lower values of the polynomial degree $k$ with respect to the approaches that use the stabilization term S$_2$. However, for the higher values of $k$, the stabilization terms S$_2$ and S$_3$ yield again similar results and very good performances are obtained when the internal DOFs \ref{DOF:2b} are employed in combination with such stabilization terms.

\begin{figure}[p]
	\centering
\vspace{-0.5\baselineskip}	\subfigure[\label{fig:ErrorHPressure_0_Q_BC3_Anisotropic_Ep_1}]
	{\includegraphics[width=.35\textwidth, height = .20\textheight]{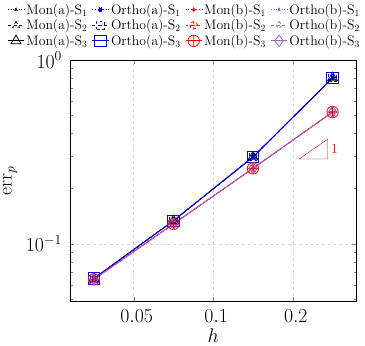}}
	\subfigure[\label{fig:ErrorHPressure_0_Q_BC3_Anisotropic_Ep_1em06}]
	{\includegraphics[width=.35\textwidth, height = .20\textheight]{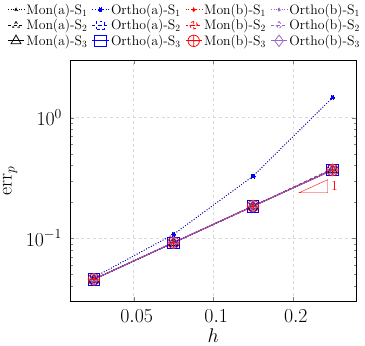}}
	
	\vspace{-0.5\baselineskip}
	\subfigure[\label{fig:ErrorHVelocity_0_Q_BC3_Anisotropic_Ep_1}]
	{\includegraphics[width=.35\textwidth, height = .20\textheight]{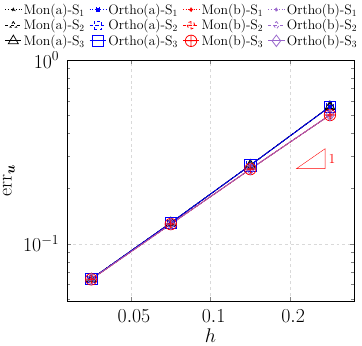}}
	\subfigure[\label{fig:ErrorHVelocity_0_Q_BC3_Anisotropic_Ep_1em06}]
	{\includegraphics[width=.35\textwidth, height = .20\textheight]{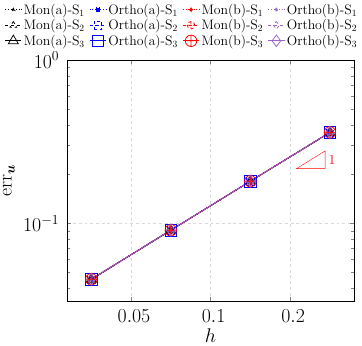}}	
	\vspace{-0.5\baselineskip}
	\caption{Test 2. Behaviour of $\mathrm{err}_{p}$ \eqref{eq:errorp} and $\mathrm{err}_{\uu}$ \eqref{eq:errorv} vs. $h$, for $k=0$ and $\Th^Q$. Left: $\epsilon = 1$.  Right: $\epsilon = 10^{-6}$. Nearly pure Neumann BCs.}
	\label{fig:errorH_Q_BC3_0_test2}
	\centering
	\subfigure[\label{fig:ErrorHPressure_0_QD_BC3_Anisotropic_Ep_1}]
	{\includegraphics[width=.35\textwidth, height = .20\textheight]{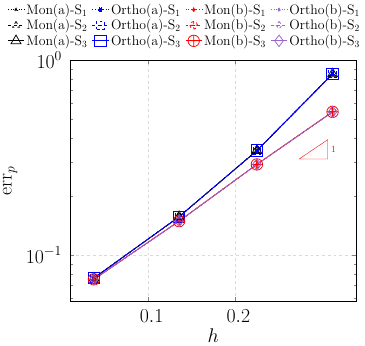}}
	\subfigure[\label{fig:ErrorHPressure_0_QD_BC3_Anisotropic_Ep_1em06}]
	{\includegraphics[width=.35\textwidth, height = .20\textheight]{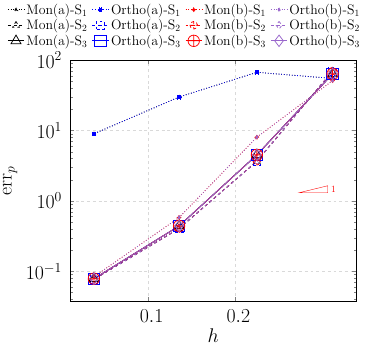}}
	
	\vspace{-0.5\baselineskip}
	\subfigure[\label{fig:ErrorHVelocity_0_QD_BC3_Anisotropic_Ep_1}]
	{\includegraphics[width=.35\textwidth, height = .20\textheight]{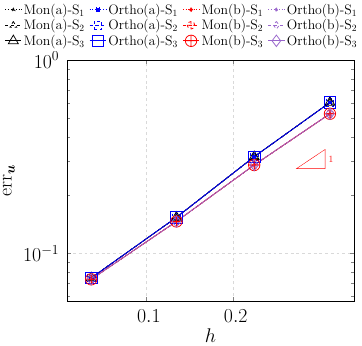}}
	\subfigure[\label{fig:ErrorHVelocity_0_QD_BC3_Anisotropic_Ep_1em06}]
	{\includegraphics[width=.35\textwidth, height = .20\textheight]{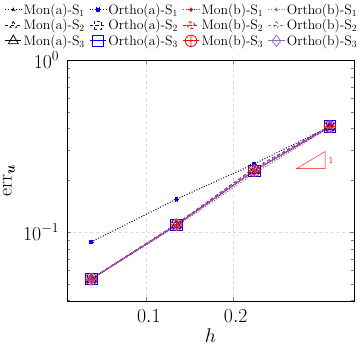}}	
	\vspace{-0.5\baselineskip}
	\caption{Test 2. Behaviour of $\mathrm{err}_{p}$ \eqref{eq:errorp} and $\mathrm{err}_{\uu}$ \eqref{eq:errorv} vs. $h$, for $k=0$ and $\Th^{DQ}$. Left: $\epsilon = 1$. Right: $\epsilon = 10^{-6}$. Nearly pure Neumann BCs.}
	\label{fig:errorH_QD_BC3_0_test2}
\end{figure}

\begin{figure}[p]
	\centering
\vspace{-0.5\baselineskip}	\subfigure[\label{fig:ErrorHPressure_5_Q_BC3_Anisotropic_Ep_1}]
	{\includegraphics[width=.35\textwidth, height = .20\textheight]{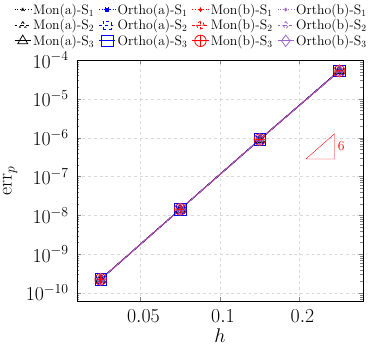}}
	\subfigure[\label{fig:ErrorHPressure_5_Q_BC3_Anisotropic_Ep_1em06}]
	{\includegraphics[width=.35\textwidth, height = .20\textheight]{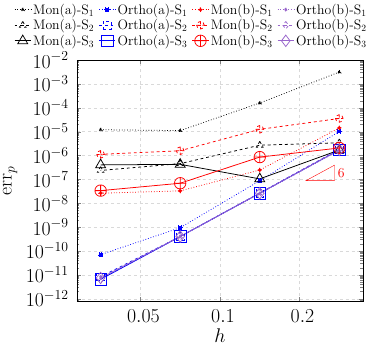}}
	
	\vspace{-0.5\baselineskip}
	\subfigure[\label{fig:ErrorHVelocity_5_Q_BC3_Anisotropic_Ep_1}]
	{\includegraphics[width=.35\textwidth, height = .20\textheight]{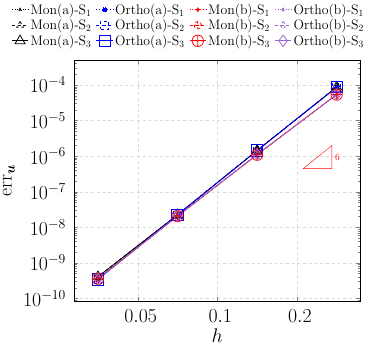}}
	\subfigure[\label{fig:ErrorHVelocity_5_Q_BC3_Anisotropic_Ep_1em06}]
	{\includegraphics[width=.35\textwidth, height = .20\textheight]{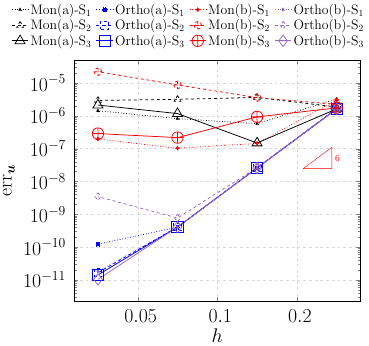}}
	\vspace{-0.5\baselineskip}
	\caption{Test 2. Behaviour of $\mathrm{err}_{p}$ \eqref{eq:errorp} and $\mathrm{err}_{\uu}$ \eqref{eq:errorv} vs. $h$, for $k=5$ and $\Th^{Q}$. Left: $\epsilon = 1$. Right: $\epsilon = 10^{-6}$. Nearly pure Neumann BCs.}
	\label{fig:errorH_Q_BC3_5_test2}

	\centering
	\subfigure[\label{fig:ErrorHPressure_5_QD_BC3_Anisotropic_Ep_1}]
	{\includegraphics[width=.35\textwidth, height = .20\textheight]{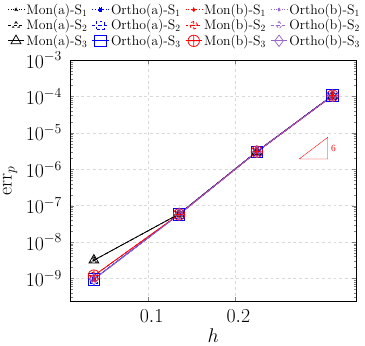}}
	\subfigure[\label{fig:ErrorHPressure_5_QD_BC3_Anisotropic_Ep_1em06}]
	{\includegraphics[width=.35\textwidth, height = .20\textheight]{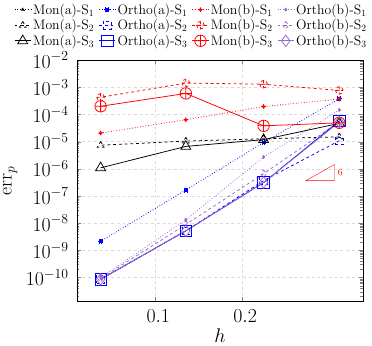}}
	
	\vspace{-0.5\baselineskip}
	\subfigure[\label{fig:ErrorHVelocity_5_QD_BC3_Anisotropic_Ep_1}]
	{\includegraphics[width=.35\textwidth, height = .20\textheight]{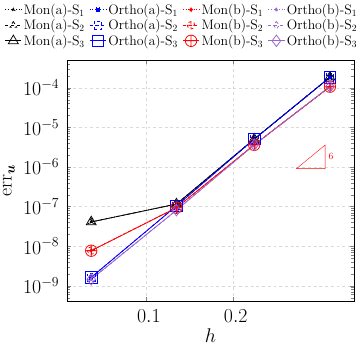}}
	\subfigure[\label{fig:ErrorHVelocity_5_QD_BC3_Anisotropic_Ep_1em06}]
	{\includegraphics[width=.35\textwidth, height = .20\textheight]{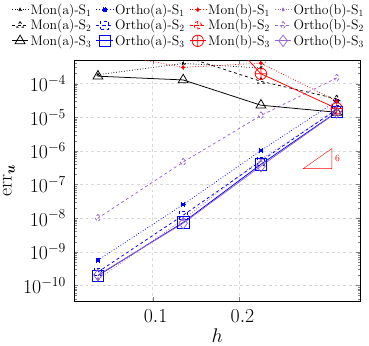}}
	\vspace{-0.5\baselineskip}
	\caption{Test 2. Behaviour of $\mathrm{err}_{p}$ \eqref{eq:errorp} and $\mathrm{err}_{\uu}$ \eqref{eq:errorv} vs. $h$, for $k=5$ and $\Th^{DQ}$. Left: $\epsilon = 1$. Right: $\epsilon = 10^{-6}$. Nearly pure Neumann BCs.}
	\label{fig:errorH_QD_BC3_5_test2}
\end{figure}

In order to analyze better such differences, in Figures \ref{fig:errorH_Q_BC3_0_test2} and \ref{fig:errorH_QD_BC3_0_test2} we report the behaviour of the errors \eqref{eq:errorp} and \eqref{eq:errorv} at decreasing values of the mesh size $h$ for the lowest polynomial degree $k=0$ and in the case of pure nearly Neumann conditions for the cartesian and the distorted quadrilateral families of meshes, respectively. In the lowest-order case, we can observe a locking phenomenon in the pressure error when distorted quadrilateral meshes are employed, as suggested by an upward shift of the error curves when $\epsilon \to 0$ and by a loss in the convergence rates, which can describe a pre-asymptotic regime \cite{locking,MANZINI2007751}. As mentioned before, the locking phenomenon is typical, generally, of the lower order methods. Indeed, looking at Figures \ref{fig:errorH_Q_BC3_5_test2} and \ref{fig:errorH_QD_BC3_5_test2}, we can note that the approaches which employed orthogonal internal DOFs show the right rates of convergence for the higher values of $k$. The monomial approaches, instead, do not converge due to ill-conditioning when $k$ is high. 

\subsection{Test 3: Two Magnetic Islands}

In the previous experiment, we considered a constant diffusion tensor with the diffusion directions aligned with the cartesian axes. Thus, the problem of anisotropy could be easily handled by choosing a cartesian family of meshes. 

\begin{figure}[]
	\centering
	\subfigure[\label{fig:campoMagnetico}]
	{\includegraphics[width=.44\textwidth, height = .22\textheight]{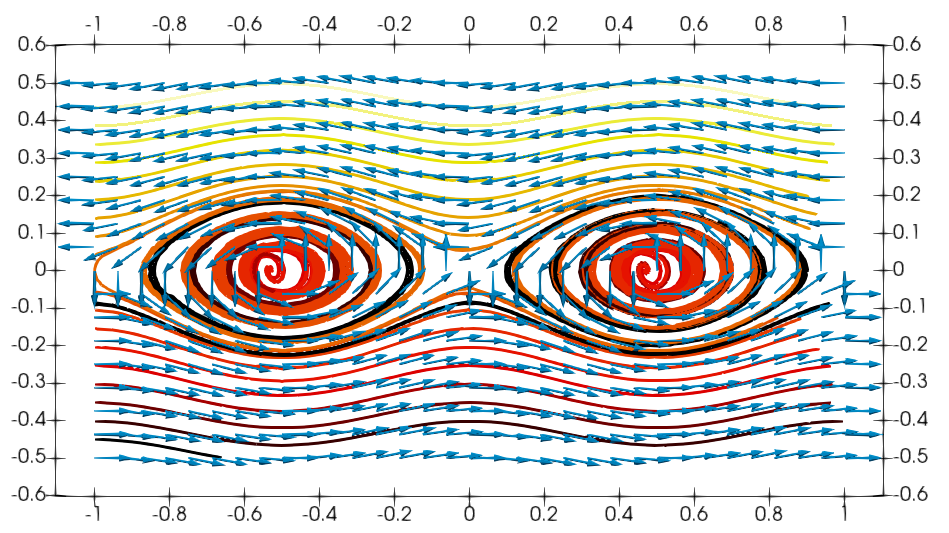}}
	\subfigure[\label{fig:sol}]
	{\includegraphics[width=.5\textwidth, height = .22\textheight]{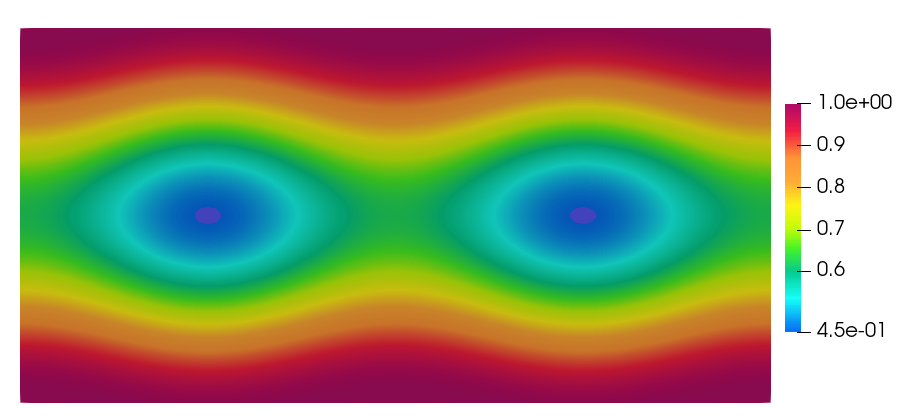}}
	\caption{Test 3. Left: the magnetic field. Right: Exact solution.}
	\label{fig:test3_solution}
\end{figure}

Now, we propose the test ``Two Magnetic Islands'' described in \cite{GREEN2022108333}, where it is almost impossible to generate an aligned mesh to solve the problem. This example models the instability phenomenon which arises in magnetized plasma for fusion applications.
More precisely, we consider a diffusion problem in $\Omega = (-1,1) \times (-0.5,0.5)$, with a diffusion tensor given by
\begin{equation}
\D(x,y) = \begin{bmatrix}
b_1(x,y) & -b_2(x,y) \\ b_2(x,y) & b_1(x,y)
\end{bmatrix} \begin{bmatrix}
D_{\vert \vert} & 0 \\ 0 & D_{\perp}
\end{bmatrix} \begin{bmatrix}
b_1(x,y) & b_2(x,y) \\ -b_2(x,y) & b_1(x,y)
\end{bmatrix},
\end{equation}
where the unit vector $\bm{b} = \begin{bmatrix}
b_1& b_2
\end{bmatrix}^T$ represents the parallel direction to the anisotropy (or to the magnetic field $\bm{B}$), while $D_{\vert \vert}$ and $D_{\perp}$ represent the parallel and the perpendicular diffusion coefficients, respectively. In this kind of application, we observe that $D_{\vert \vert}$ can be greater than $D_{\perp}$ by a factor of $10^{12}$ \cite{GREEN2022108333}. Let us now define the equilibrium magnetic field
\begin{equation}
\bm{B}(x,y) = \begin{bmatrix}
- \pi \sin(\pi y)\\
\frac{2 \pi}{10} \sin\left(2 \pi \left(x - \frac{3}{2}\right)\right)
\end{bmatrix}
\end{equation}  
which is shown in Figure \ref{fig:campoMagnetico}. By looking at this figure, we note that the magnetic field results to be the zero-vector in the center of the ``magnetic islands'' (the $\mathcal{O}$-points) and where the field lines cross each other (the so-called $\mathcal{X}$-points). In all the other points, we can define
\begin{equation}
\bb(x,y) = \frac{\bm{B}(x,y)}{\| \bm{B}(x,y)\|}
\end{equation}
and compute
\begin{equation}
\D(x,y)^{-1} = \begin{bmatrix}
b_1(x,y) & -b_2(x,y) \\ b_2(x,y) & b_1(x,y)
\end{bmatrix} \begin{bmatrix}
\frac{1}{D_{\vert \vert}} & 0 \\ 0 & \frac{1}{D_{\perp}}
\end{bmatrix} \begin{bmatrix}
b_1(x,y) & b_2(x,y) \\ -b_2(x,y) & b_1(x,y)
\end{bmatrix}.
\end{equation}
We further fix $D_{\perp} = 1$, while $D_{\vert \vert} \in \{1,10^{4},10^{8}\}$. 

We evaluate the performances of the aforementioned approaches on a family $\Th^S = \{\Th[i]^S\}_{i=1}^{4}$ of four squared meshes, which are characterized by an edge length decreasing as $\frac{1}{2^{i+1}}$, with $i=1,\dots,4$. We note that both the $\mathcal{O}$-points and $\mathcal{X}$-points represent vertices of the tessellation in each refinement.
Furthermore, we define the forcing term and the boundary conditions in such a way the exact solution is 
\begin{equation}
p(x,y) = \cos\left(\frac{1}{10}\cos\left(2 \pi \left(x - \frac{3}{2}\right)\right)+\cos(\pi y)\right),
\end{equation}
which is shown in Figure \ref{fig:sol}. We test two cases, characterized by different boundary conditions, namely
\begin{itemize}
\item pure Dirichlet boundary conditions $\Gamma_D = \Gamma$;
\item mixed boundary conditions, with
\begin{linenomath}
\begin{equation*}
\Gamma_N = \{(x,y): x = -1  \text{ or } x = 1\}.
\end{equation*}
\end{linenomath}
\end{itemize}
We note that the velocity field does not depend on the parameter $D_{\vert \vert}$.

In this experiment, we test three possible choices for the stabilization term, namely
\begin{itemize}
\item S$_1$: the dofi-dofi stabilization with $C_{\D^{-1}} = \frac{1}{D_{\vert \vert}}$. 
\item S$_2$: the D-recipe stabilization with $C_{\D^{-1}} = 1$.
\item S$_3$: a D-recipe stabilization term with
\begin{linenomath}
\begin{equation*}
S_{ii} = \vert E \vert \begin{cases}
\max(\nn_{e_i} \cdot \D^{-1}(\x_{e_i}) \nn_{e_i}, (\KK_C^E)_{ii})& \text{if } i \text{ is a boundary DOF}\\
0& \text{if } i \text{ is an internal DOF}
\end{cases},
\end{equation*}
\end{linenomath}
where $\x_{e_i}$ and $\nn_{e_i}$ are the midpoint and the unit outward normal vector to the edge $e_i$ related to the boundary DOF $i$. This stabilization term, inspired by \cite{GIORGIANI2020107375}, aims to take into account the actual strength of the normal contribution of the parallel diffusion on each edge.
\end{itemize}

\begin{figure}[p]
	\centering
	\subfigure[\label{fig:ErrorL2Pressure_Q16x8_Anisotropic_Ep_1}]
	{\includegraphics[width=.32\textwidth, height = .19\textheight]{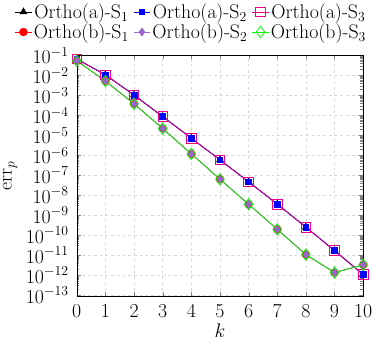}}
	\subfigure[\label{fig:ErrorL2Pressure_Q16x8_Anisotropic_Ep_10000}]
	{\includegraphics[width=.32\textwidth, height = .19\textheight]{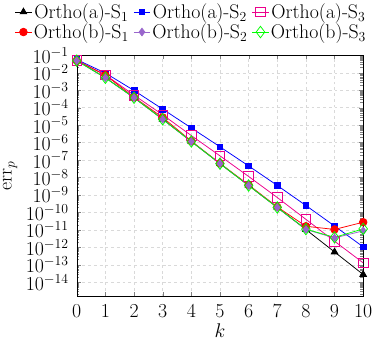}}
	\subfigure[\label{fig:ErrorL2Pressure_Q16x8_Anisotropic_Ep_100000000}]
	{\includegraphics[width=.32\textwidth, height = .19\textheight]{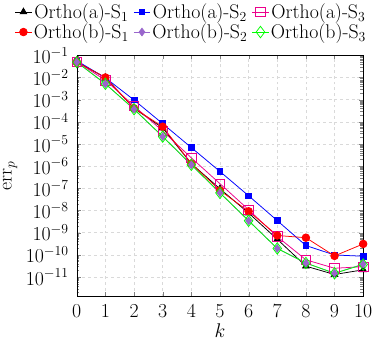}}
	
	\subfigure[\label{fig:ErrorL2Velocity_Q16x8_Anisotropic_Ep_1}]
	{\includegraphics[width=.32\textwidth, height = .19\textheight]{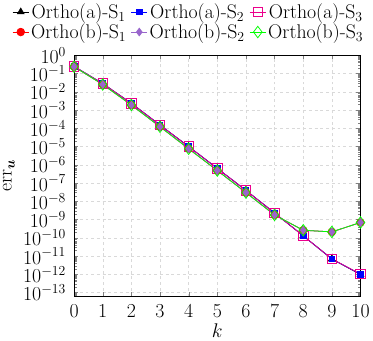}}
	\subfigure[\label{fig:ErrorL2Velocity_Q16x8_Anisotropic_Ep_10000}]
	{\includegraphics[width=.32\textwidth, height = .19\textheight]{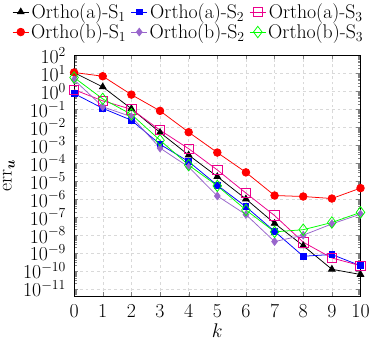}}
	\subfigure[\label{fig:ErrorL2Velocity_Q16x8_Anisotropic_Ep_100000000}]
	{\includegraphics[width=.32\textwidth, height = .19\textheight]{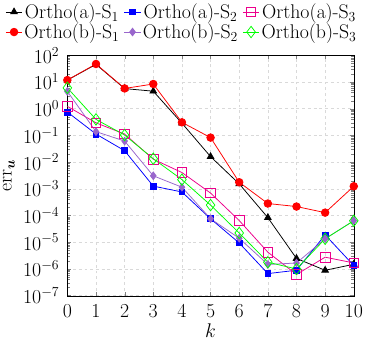}}
	\caption{Test 3. Behaviour of $\mathrm{err}_p$ \eqref{eq:errorp} and $\mathrm{err}_{\bm{u}}$ \eqref{eq:errorv} vs. $k$, for $\Th[2]^S$. Left: $D_{\vert \vert} = 1$. Center: $D_{\vert \vert}= 10^{4}$. Right: $D_{\vert \vert} = 10^{8}$. Dirichlet BCs.}
	\label{fig:error_test3_dirichlet}

	\centering
	\subfigure[\label{fig:Mixed_ErrorL2Pressure_Q16x8_Anisotropic_Ep_1}]
	{\includegraphics[width=.32\textwidth, height = .19\textheight]{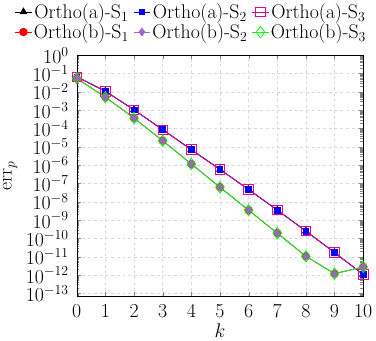}}
	\subfigure[\label{fig:Mixed_ErrorL2Pressure_Q16x8_Anisotropic_Ep_10000}]
	{\includegraphics[width=.32\textwidth, height = .19\textheight]{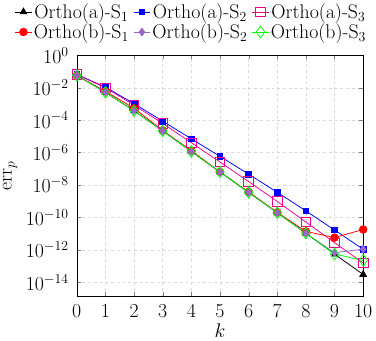}}
	\subfigure[\label{fig:Mixed_ErrorL2Pressure_Q16x8_Anisotropic_Ep_100000000}]
	{\includegraphics[width=.32\textwidth, height = .19\textheight]{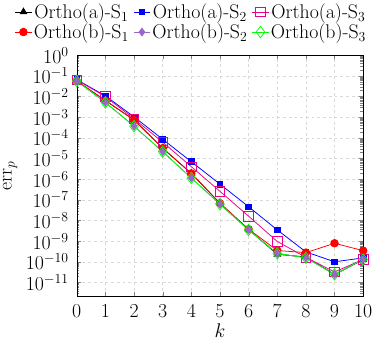}}
	
	\subfigure[\label{fig:Mixed_ErrorL2Velocity_Q16x8_Anisotropic_Ep_1}]
	{\includegraphics[width=.32\textwidth, height = .19\textheight]{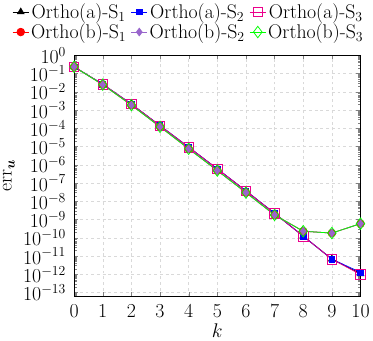}}
	\subfigure[\label{fig:Mixed_ErrorL2Velocity_Q16x8_Anisotropic_Ep_10000}]
	{\includegraphics[width=.32\textwidth, height = .19\textheight]{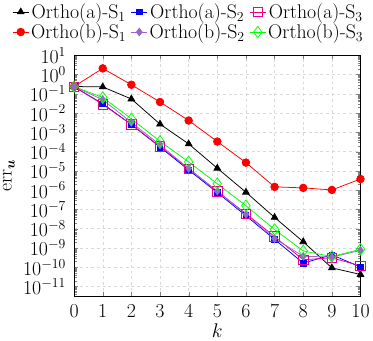}}
	\subfigure[\label{fig:Mixed_ErrorL2Velocity_Q16x8_Anisotropic_Ep_100000000}]
	{\includegraphics[width=.32\textwidth, height = .19\textheight]{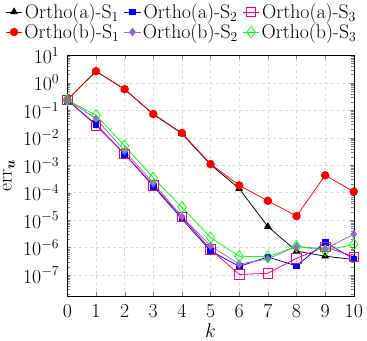}}
	\caption{Test 3. Behaviour of $\mathrm{err}_p$ \eqref{eq:errorp} and $\mathrm{err}_{\bm{u}}$ \eqref{eq:errorv} vs. $k$, for $\Th[2]^S$. Left: $D_{\vert \vert} = 1$. Center: $D_{\vert \vert}= 10^{4}$. Right: $D_{\vert \vert} = 10^{8}$. Mixed BCs.}
	\label{fig:error_test3_mixed}
\end{figure}

Figures \ref{fig:error_test3_dirichlet} and \ref{fig:error_test3_mixed} show the behaviour of the errors \eqref{eq:errorp} and \eqref{eq:errorv} at varying of the polynomial degree $k$ for the second refinement $\Th[2]^S$ when the Dirichlet and mixed boundary conditions are imposed. We decide to report only the behaviour of the Ortho approaches in these figures in order to try to better highlight differences between the employment of boundary DOFs \ref{DOF:1a} and \ref{DOF:1b}.

\begin{figure}[p]
	\centering
	\subfigure[\label{fig:ErrorHPressure_0_Q_Anisotropic_Ep_1}]
	{\includegraphics[width=.32\textwidth, height = .19\textheight]{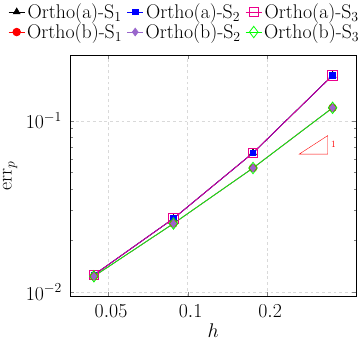}}
	\subfigure[\label{fig:ErrorHPressure_0_Q_Anisotropic_Ep_10000}]
	{\includegraphics[width=.32\textwidth, height = .19\textheight]{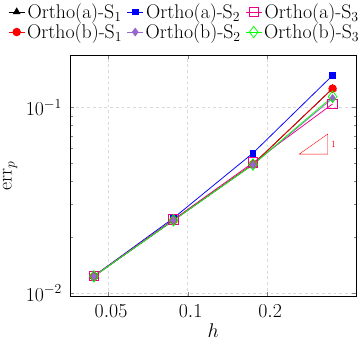}}
	\subfigure[\label{fig:ErrorHPressure_0_Q_Anisotropic_Ep_100000000}]
	{\includegraphics[width=.32\textwidth, height = .19\textheight]{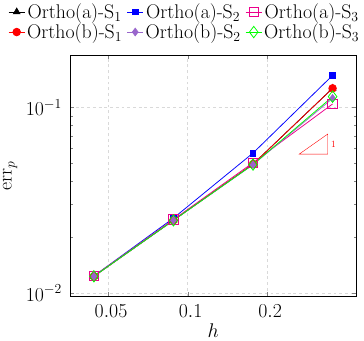}}
	
	\subfigure[\label{fig:ErrorHVelocity_0_Q_Anisotropic_Ep_1}]
	{\includegraphics[width=.32\textwidth, height = .19\textheight]{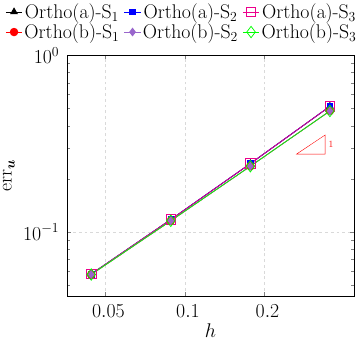}}
	\subfigure[\label{fig:ErrorHVelocity_0_Q_Anisotropic_Ep_10000}]
	{\includegraphics[width=.32\textwidth, height = .19\textheight]{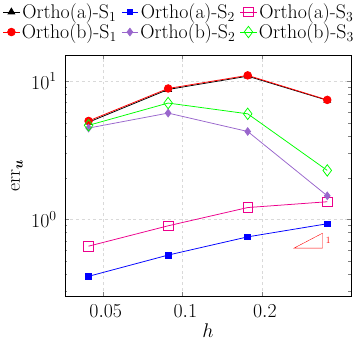}}
	\subfigure[\label{fig:ErrorHVelocity_0_Q_Anisotropic_Ep_100000000}]
	{\includegraphics[width=.32\textwidth, height = .19\textheight]{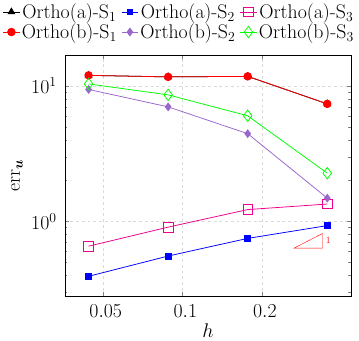}}
	
	\caption{Test 3. Behaviour of $\mathrm{err}_{p}$ \eqref{eq:errorp} and $\mathrm{err}_{\uu}$ \eqref{eq:errorv} vs. $h$, for $k = 0$. Left: $D_{\vert \vert} = 1$. Center: $D_{\vert \vert}= 10^{4}$. Right: $D_{\vert \vert} = 10^{8}$. Dirichlet BCs.}
	\label{fig:errorH_Q_Dirichlet_0_test3}

	\centering
	\subfigure[\label{fig:Mixed_ErrorHPressure_0_Q_Anisotropic_Ep_1}]
	{\includegraphics[width=.32\textwidth, height = .19\textheight]{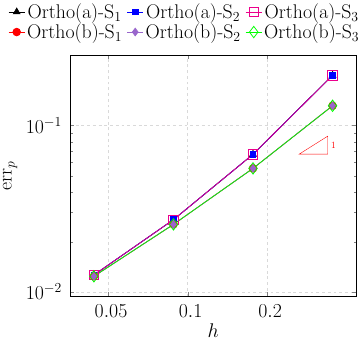}}
	\subfigure[\label{fig:Mixed_ErrorHPressure_0_Q_Anisotropic_Ep_10000}]
	{\includegraphics[width=.32\textwidth, height = .19\textheight]{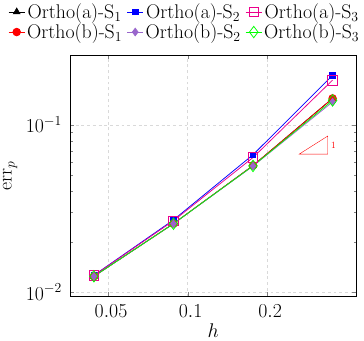}}
	\subfigure[\label{fig:Mixed_ErrorHPressure_0_Q_Anisotropic_Ep_100000000}]
	{\includegraphics[width=.32\textwidth, height = .19\textheight]{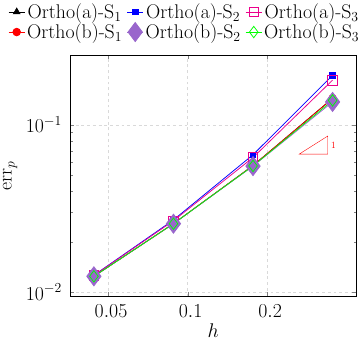}}
	
	\subfigure[\label{fig:Mixed_ErrorHVelocity_0_Q_Anisotropic_Ep_1}]
	{\includegraphics[width=.32\textwidth, height = .19\textheight]{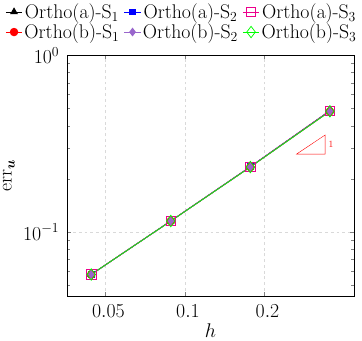}}
	\subfigure[\label{fig:Mixed_ErrorHVelocity_0_Q_Anisotropic_Ep_10000}]
	{\includegraphics[width=.32\textwidth, height = .19\textheight]{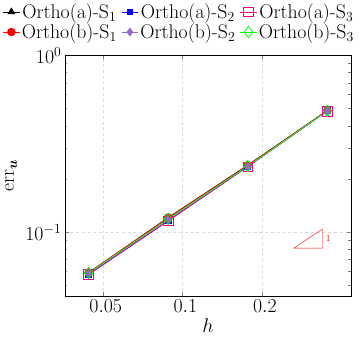}}
	\subfigure[\label{fig:Mixed_ErrorHVelocity_0_Q_Anisotropic_Ep_100000000}]
	{\includegraphics[width=.32\textwidth, height = .19\textheight]{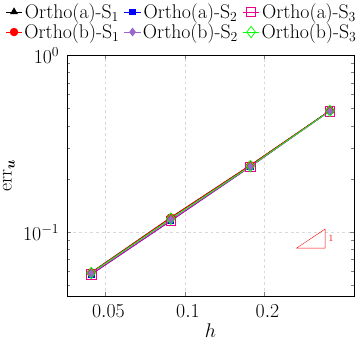}}
	
	\caption{Test 3. Behaviour of $\mathrm{err}_{p}$ \eqref{eq:errorp} and $\mathrm{err}_{\uu}$ \eqref{eq:errorv} vs. $h$, for $k = 0$. Left: $D_{\vert \vert} = 1$. Center: $D_{\vert \vert}= 10^{4}$. Right: $D_{\vert \vert} = 10^{8}$. Mixed BCs.}
	\label{fig:errorH_Q_Mixed_0_test3}
\end{figure}

\begin{figure}[p]
	\centering
	\subfigure[\label{fig:ErrorHPressure_2_Q_Anisotropic_Ep_1}]
	{\includegraphics[width=.32\textwidth, height = .19\textheight]{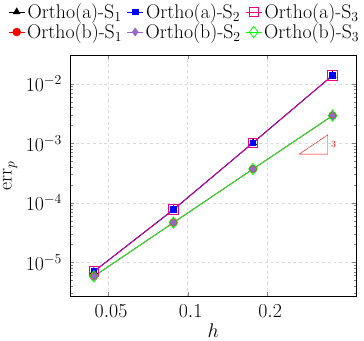}}
	\subfigure[\label{fig:ErrorHPressure_2_Q_Anisotropic_Ep_10000}]
	{\includegraphics[width=.32\textwidth, height = .19\textheight]{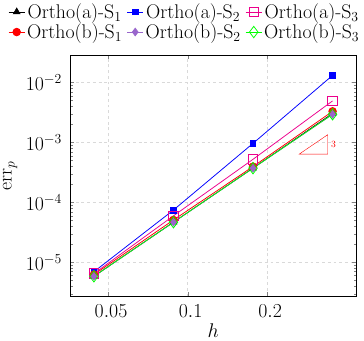}}
	\subfigure[\label{fig:ErrorHPressure_2_Q_Anisotropic_Ep_100000000}]
	{\includegraphics[width=.32\textwidth, height = .19\textheight]{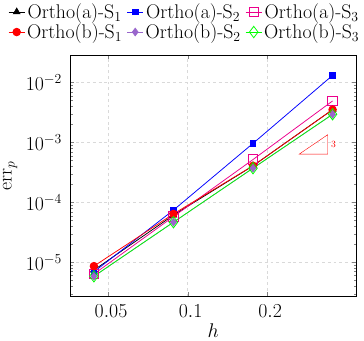}}
	
	\subfigure[\label{fig:ErrorHVelocity_2_Q_Anisotropic_Ep_1}]
	{\includegraphics[width=.32\textwidth, height = .19\textheight]{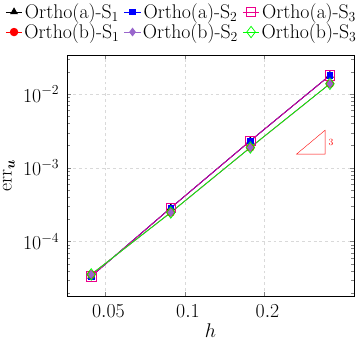}}
	\subfigure[\label{fig:ErrorHVelocity_2_Q_Anisotropic_Ep_10000}]
	{\includegraphics[width=.32\textwidth, height = .19\textheight]{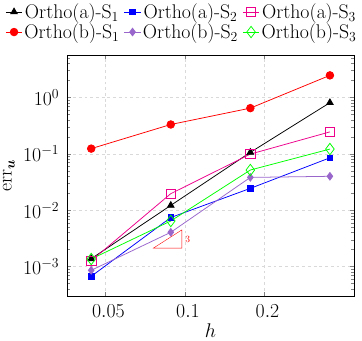}}
	\subfigure[\label{fig:ErrorHVelocity_2_Q_Anisotropic_Ep_100000000}]
	{\includegraphics[width=.32\textwidth, height = .19\textheight]{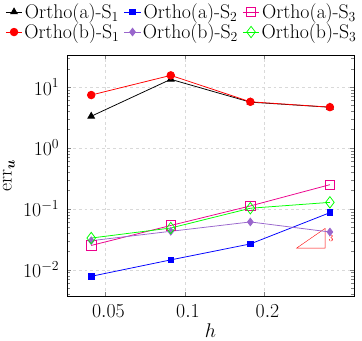}}
	
	\caption{Test 3. Behaviour of the $\mathrm{err}_{p}$ \eqref{eq:errorp} and $\mathrm{err}_{\uu}$ \eqref{eq:errorv} vs. $h$, for $k = 2$. Left: $D_{\vert \vert} = 1$. Center: $D_{\vert \vert}= 10^{4}$. Right: $D_{\vert \vert} = 10^{8}$. Dirichlet BCs.}
	\label{fig:errorH_Q_Dirichlet_2_test3}

	\centering
	\subfigure[\label{fig:Mixed_ErrorHPressure_2_Q_Anisotropic_Ep_1}]
	{\includegraphics[width=.32\textwidth, height = .19\textheight]{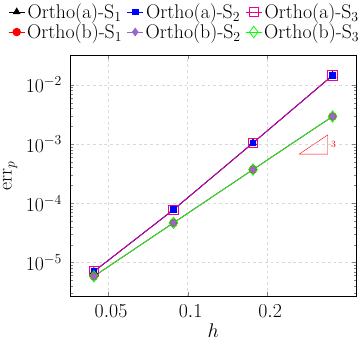}}
	\subfigure[\label{fig:Mixed_ErrorHPressure_2_Q_Anisotropic_Ep_10000}]
	{\includegraphics[width=.32\textwidth, height = .19\textheight]{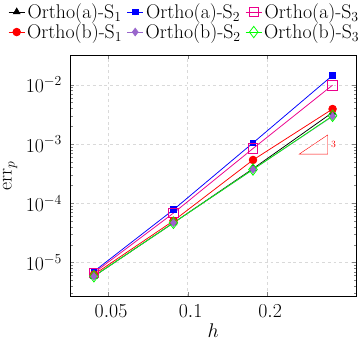}}
	\subfigure[\label{fig:Mixed_ErrorHPressure_2_Q_Anisotropic_Ep_100000000}]
	{\includegraphics[width=.32\textwidth, height = .19\textheight]{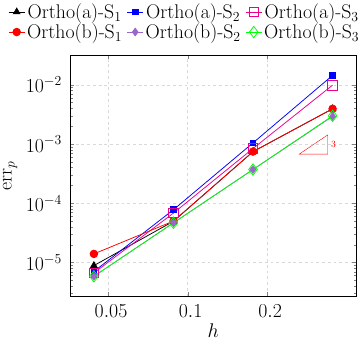}}
	
	\subfigure[\label{fig:Mixed_ErrorHVelocity_2_Q_Anisotropic_Ep_1}]
	{\includegraphics[width=.32\textwidth, height = .19\textheight]{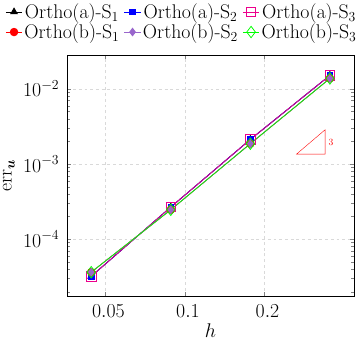}}
	\subfigure[\label{fig:Mixed_ErrorHVelocity_2_Q_Anisotropic_Ep_10000}]
	{\includegraphics[width=.32\textwidth, height = .19\textheight]{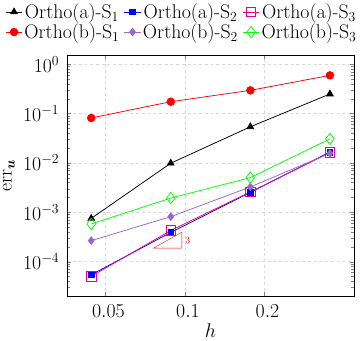}}
	\subfigure[\label{fig:Mixed_ErrorHVelocity_2_Q_Anisotropic_Ep_100000000}]
	{\includegraphics[width=.32\textwidth, height = .19\textheight]{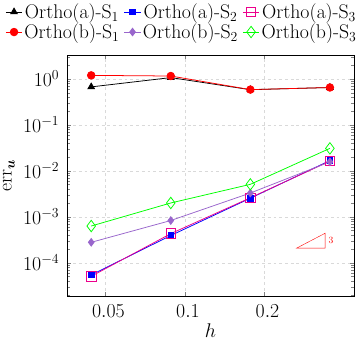}}
	
	\caption{Test 3. Behaviour of $\mathrm{err}_{p}$ \eqref{eq:errorp} and $\mathrm{err}_{\uu}$ \eqref{eq:errorv} vs. $h$, for $k = 2$. Left: $D_{\vert \vert} = 1$. Center: $D_{\vert \vert}= 10^{4}$. Right: $D_{\vert \vert} = 10^{8}$. Mixed BCs.}
	\label{fig:errorH_Q_Mixed_2_test3}
\end{figure}

We observe that all approaches show the right behaviour in terms of the relative pressure error \eqref{eq:errorp} in the case of both Dirichlet and mixed boundary conditions. This appears also evident when observing the behaviour of the pressure error in terms of $h$ in Figures \ref{fig:errorH_Q_Dirichlet_0_test3}, \ref{fig:errorH_Q_Mixed_0_test3} \ref{fig:errorH_Q_Dirichlet_2_test3} and \ref{fig:errorH_Q_Mixed_2_test3} for the lowest order $k=0$ and for the polynomial degree $k =2$. From these figures we can note that, as usual, the choice of boundary DOFs \ref{DOF:1b} is characterized by smaller pressure error constants with respect to the choice  \ref{DOF:1a}. Furthermore, approaches that employ \ref{DOF:1b} seem to be less sensitive to the choice of the stabilization term 
than approaches which exploit boundary DOFs \ref{DOF:1a}.

However, the same conclusions do not hold true when dealing with the relative velocity error \eqref{eq:errorv}. Indeed, we first can note that switching off the stabilization by choosing the stabilization term S$_1$ when $D_{\vert \vert}$ is big enough generally does not lead to good results in terms of the velocity error. Furthermore, we note that in order to achieve good results in terms of the velocity error, it is very important to enforce the velocity on the boundary by imposing strong Neumann boundary conditions when high values of $D_{\vert \vert}$ are considered. In this way, it is possible to obtain the right convergence rates in terms of the mesh size of both the pressure and the velocity errors as can be seen in Figures \ref{fig:errorH_Q_Dirichlet_0_test3} and \ref{fig:errorH_Q_Mixed_0_test3}.

Finally, we observe that, in this test case, the Ortho (a) approach seems to perform better than the Ortho (b) approach in terms of velocity error when highly anisotropic cases are taken into account.
  

\section{Conclusions}

In this paper, we carried out the analysis of the robustness of the mixed Virtual Element Method when problems characterized by highly anisotropic diffusion tensors are considered. Furthermore, a new set of boundary degrees of freedom based on moments computed against an $L^2([0,1])$-orthonormal basis is also introduced.

Here, we report the results obtained on a set of benchmark problems by resorting to various approaches which differ for the sets of both the internal and the boundary degrees of freedom. For each benchmark problem, we propose different kinds of the stabilization term and we test the sensitivity of each proposed approach to the choice of the stabilization term in terms of both the condition number of the system matrix and of the errors \eqref{eq:errorp} and \eqref{eq:errorv}. 

In particular, the new set of boundary degrees of freedom seems to be more favourable in terms of errors by leading to a downward shift of the error curves, although, this choice generally does not ensure obtaining an improvement in the conditioning of $\KK$. Indeed, the condition number of the system matrix seems to be mainly controlled by the choice of internal DOFs and by the anisotropic ratio.

Finally, the D-recipe version of the stabilization term with unit constant seems to be a good alternative to build a robust method for highly anisotropic diffusion problems.

\section*{Acknowledgments}

The author S.B. kindly acknowledges partial financial support provided by PRIN project “Advanced polyhedral discretisations of heterogeneous PDEs for multiphysics problems” (No. 20204LN5N5\_003) and by PNRR M4C2 project of CN00000013 National Centre for HPC, Big Data and Quantum Computing (HPC) (CUP: E13C22000990001). 
The author S.S. kindly acknowledges partial financial support provided by INdAM-GNCS through project “Sviluppo ed analisi di Metodi agli Elementi Virtuali per processi accoppiati su geometrie complesse” and that this  publication is part of the project NODES which has received funding from the MUR-M4C2 1.5 of PNRR with grant agreement no. ECS00000036. 
The author G.T. kindly acknowledges financial support provided by the MIUR programme ``Programma Operativo Nazionale Ricerca e Innovazione 2014 - 2020'' 
~ (CUP: E11B21006490005). Computational resources are partially supported by SmartData@polito. The authors are members of the Italian INdAM-GNCS research group.


\bibliographystyle{IEEEtran}
\bibliography{biblio.bib}

\end{document}